\documentclass{scrartcl}
\usepackage{enumitem}
\usepackage{mathtools}

\usepackage{tikz}
\usetikzlibrary{decorations.pathreplacing}

\usepackage{booktabs} 
\usepackage{amsmath}
\usepackage{amsfonts}
\usepackage{amssymb}
\usepackage{amsthm}
\usepackage{dsfont}
\usepackage[T1]{fontenc}
\usepackage[latin1]{inputenc}
\usepackage{epsfig}
\usepackage{graphicx}
\usepackage{bm}
\usepackage[english]{babel}

\usepackage[boxed, lined]{algorithm2e}
\usepackage{float}
\usepackage{comment}

\newcommand{\bj}{\mathbf{j}}

\newcommand{\ba}{\mathbf{a}}

\newcommand{\bx}{\mathbf{x}}

\newcommand{\D}{{\mathcal{D}}}

\newcommand{\F}{{\mathcal{F}}}
\newcommand{\G}{{\mathcal{G}}}
\newcommand{\HH}{{\mathcal{H}}}

\newcommand{\Nu}{{\mathcal{N}}}

\newcommand{\N}{\mathbb{N}}

\newcommand{\R}{\mathbb{R}}
\newcommand{\Z}{\mathbb{Z}}
\newcommand{\Rd}{\mathbb{R}^d}

\newcommand{\beq}{\begin{eqnarray*}}

\newcommand{\eeq}{\end{eqnarray*}}

\newcommand{\beqm}{\begin{eqnarray}}

\newcommand{\eeqm}{\end{eqnarray}}

\newtheorem{theorem}{Theorem}

\newtheorem{lemma}{Lemma}
\newtheorem{definition}{Definition}

\DeclareMathOperator*{\argmax}{arg\,max}

\DeclareMathOperator{\sgn}{sgn}

\newcommand{\EXP}{{\mathbf E}}
\newcommand{\PROB}{{\mathbf P}}

\renewcommand{\P}{{\cal P}}
\renewcommand{\bf}{\normalfont \bfseries}
\renewcommand{\it}{\normalfont \itshape}

\newcommand{\secrev}[2]{{#2}}
\newcommand{\rev}[2]{{#2}}
\newcommand{\gozde}[1]{#1}

\allowdisplaybreaks
\begin{document}
\begin{center}

  {\LARGE \bf
    On the rate of convergence of a classifier based on
    a Transformer encoder}\footnote{
Running title: {\it Rate of convergence of a Transformer classifier}}
\vspace{0.5cm}

Iryna Gurevych$^\dagger$, Michael Kohler$^{*,}\footnote{Corresponding author. Tel:
  +49-6151-16-23382}$, G\"{o}zde G\"{u}l \c{S}ahin$^\dagger$ \\

{\it 
  $^*$ Department of Mathematics, Technical University of Darmstadt,
Schlossgartenstr. 7, 64289 Darmstadt, Germany,
email: kohler@mathematik.tu-darmstadt.de \\
$^\dagger$ Ubiquitous Knowledge Processing (UKP) Lab, Technical University of Darmstadt, Germany,
email: gurevych@ukp.informatik.tu-darmstadt.de,
sahin@ukp.informatik.tu-darmstadt.de.
}

\end{center}
\vspace{0.5cm}

\begin{center}
November 22, 2021
\end{center}
\vspace{0.5cm}

\noindent
    {\bf Abstract}\\
    Pattern recognition based on a high-dimensional
    predictor is considered. A classifier is defined
    which is based on a Transformer encoder.
    The rate of convergence of the misclassification
    probability of the classifier towards the
    optimal misclassification probability is analyzed.
    It is shown that this classifier is able to
    circumvent the curse of dimensionality provided
    the aposteriori probability satisfies a suitable
    hierarchical composition model. Furthermore,
    the difference between Transformer classifiers
    analyzed theoretically in this paper and Transformer
    classifiers used nowadays in practice are illustrated
    by considering classification problems in natural
    language processing.
    \vspace*{0.2cm}


\vspace*{0.2cm}


\section{Introduction}
\label{se1}
Deep learning has achieved impressive progress in natural language
processing (NLP), e.g. in the areas of understanding,
summarizing or generating text, see., e.g., Young et al. (2018)
and the literature cited therein.
Among the most successful techniques in this field are Transformers
introduced by Vaswani et al. (2017). They used a decoder-encoder
structure based on multi-head attention and piecewise feedforward layers
and achieved a top performance in applications in machine translation.
Devlin et al. (2019) describe how Transformers can be combined
with unsupervised pre-training such that the same pre-trained
Transformer encoder can be fine-tuned to a variety of
natural language processing tasks.

Besides the huge practical success of these estimates, their theoretical
performance has not been studied intensively until now. This is in sharp
contrast to deep neural networks, where various results concerning
the approximation power of deep neural networks (cf., e.g., 
Yarotsky (2017),  Yarotsky and Zhevnerchuck (2020),
Lu et al.  (2020), Langer (2021b)
and the literature cited 
therein)
or concerning the statistical risk of
corresponding estimates
(cf., e.g.,
Bauer and Kohler (2019), Kohler and Krzy\.zak (2017),
Schmidt-Hieber (2019), 
Kohler and Langer (2020), 
Langer (2021a), Imaizumi and Fukumizu (2019), Suzuki (2018), Suzuki 
and Nitanda (2019),  and the literature cited 
therein)
have been derived in the last few years.

In this paper we aim to shed light on the theoretical performance
of Transformers. To do this, we focus on pattern recognition and
consider estimates based on a Transformer encoder. Here we simplify
the learning problem by defining our estimate as
a plug-in classification rule based on
an abstract least
squares estimates, i.e., we ignore the gradient descent usually applied
in practice as this gradient descent is even for deep neural networks
nowadays not well understood.
After this simplification the main remaining challenges in view of
a theoretical understanding are the
approximation properties and the generalization abilities of Transformer
encoders.

We study these estimates in the context of pattern recognition.
Here, $(X,Y)$, $(X_1,Y_1)$, \dots, $(X_n,Y_n)$
are independent and identically distributed random variables
with values in $\R^{d \cdot l} \times \{0,1\}$, and given the data set
\[
\D_n = \{ (X_1,Y_1), \dots, (X_n,Y_n)\}
\]
the goal is to construct a classifier
\[
\eta_n(\cdot) = \eta_n(\cdot,\D_n):\R^{d \cdot l} \rightarrow \{0,1\}
\]
such that the misclassification probability
\[
\PROB \{ \eta_n(X) \neq Y | \D_n\}
\]
is as small as possible.
Here the predictor variable $X$ describes the encoding of a sequence
of lenght $l$ consisting of words or tokens, and each word or token
is encoded by a value in $\R^d$.

Let
\begin{equation}
\label{se1eq1}
m(x) = \PROB\{ Y=1|X=x\}
\quad
( x \in \R^{d \cdot l})
\end{equation}
be the so--called aposteriori probability of class 1. Then
\[
\eta^*(x)=
\begin{cases}
  1, & \mbox{if } m(x) > \frac{1}{2} \\
  0, & \mbox{elsewhere}
  \end{cases}
\]
is the so--called Bayes classifier, i.e., it satisfies
\[
\PROB\{ \eta^*(X) \neq Y\}
=
\min_{\eta : \R^{d \cdot l} \rightarrow \{0,1\}}
\PROB\{\eta(X) \neq Y\}
\]
(cf., e.g., Theorem 2.1 in Devroye, Gy\"orfi and Lugosi (1996)).

In the sequel we try to derive upper bounds on
\begin{eqnarray}
  \label{se1eq2}
  &&
  \EXP \left\{
\PROB \{ \eta_n(X) \neq Y | \D_n\}
-
\PROB\{\eta^*(X) \neq Y\}
\right\}
\nonumber \\
&&
=
\PROB \{ \eta_n(X) \neq Y \}
-
\min_{\eta : \R^{d \cdot l} \rightarrow \{0,1\}}
\PROB\{\eta(X) \neq Y\}
.
\end{eqnarray}
It is well-known that in order to derive
nontrivial rate of convergence results
on the difference between the misclassification risk of any
estimate and the minimal possible value it is necessary
to restrict the class of distributions (cf.,
Cover (1968) and Devroye (1982)). In this context we will
assume that the aposteriori probability is smooth, and for this
we will use our next definition.

\begin{definition}
\label{intde2} 
  Let $p=q+s$ for some $q \in \N_0$ and $0< s \leq 1$.
A function $m:\R^{d \cdot l} \rightarrow \R$ is called
$(p,C)$-smooth, if for every $\bm{\alpha}=(\alpha_1, \dots, \alpha_{d \cdot l}) \in
\N_0^{d \cdot l}$
with $\sum_{j=1}^{d \cdot l} \alpha_j = q$ the partial derivative
$\partial^q m/(\partial x_1^{\alpha_1}
\dots
\partial x_{d \cdot l}^{\alpha_{d \cdot l}}
)$
exists and satisfies
\[
\left|
\frac{
\partial^q m
}{
\partial x_1^{\alpha_1}
\dots
\partial x_{d \cdot l}^{\alpha_{d \cdot l}}
}
(x)
-
\frac{
\partial^q m
}{
\partial x_1^{\alpha_1}
\dots
\partial x_{d \cdot l}^{\alpha_{d \cdot l}}
}
(z)
\right|
\leq
C
\|\bold{x}-\bold{z}\|^s
\]
for all $\bold{x},\bold{z} \in \R^{d \cdot l}$, where $\Vert\cdot\Vert$ denotes the Euclidean norm.
\end{definition}

In order to be able to show good rate of convergences
even for high-dimensional predictors we use a hierarchical composition
model as in Schmidt-Hieber (2019), where the aposteriori probability
is represented by a composition of several functions and where each
of these functions depends only on a few variables. We use the
following definition of Kohler and Langer (2020) to formalize
this assumption.

\begin{definition}
\label{de2}
Let $d,l \in \N$, $m: \R^{d \cdot l} \to \R$ and let
$\P$ be a subset
of $(0,\infty) \times \N$.\\
\noindent
\textbf{a)}
We say that $m$ satisfies a hierarchical composition model of level $0$
with order and smoothness constraint $\mathcal{P}$, if there exists $K \in \{1, \dots, {d \cdot l}\}$ such that
\[
m(\bold{x}) = x^{(K)} \quad \mbox{for all } \bold{x}
= (x^{(1)}, \dots, x^{({d \cdot l})})^{\top} \in \R^{d \cdot l}.
\]
\noindent
\textbf{b)}
Let $\kappa \in \N_0$.
We say that $m$ satisfies a hierarchical composition model
of level $\kappa+1$ with order and smoothness constraint $\mathcal{P}$, if there exist $(p,K)  \in \P$, $C>0$, \linebreak $g: \R^{K} \to \R$ and $f_{1}, \dots, f_{K}: \Rd \to \R$, such that
$g$ is $(p,C)$--smooth,
$f_{1}, \dots, f_{K}$ satisfy a  hierarchical composition model of level $\kappa$
with order and smoothness constraint $\mathcal{P}$
and 
\[m(\bold{x})=g(f_{1}(\bold{x}), \dots, f_{K}(\bold{x})) \quad \mbox{for all } \bold{x} \in \R^{d \cdot l}.\]
\end{definition}
\noindent
Let $\HH(\kappa,\P)$ be the set of all functions $m:\R^{d \cdot l} \rightarrow \R$
which satisfy a  hierarchical composition model of level $\kappa$
with order and smoothness constraint $\mathcal{P}$.

It was shown in Bauer and Kohler (2019), 
Schmidt-Hieber (2019) and 
Kohler and Langer (2020) that deep neural networks are able to
circumvent the curse of dimensionality in case that the function
to be estimated is
contained
in a suitably defined hierarchical composition model.
The main contribution of this paper is to show that classifiers
based on Transformer encoders, which we will introduce in the next
section, have this property, too. More precisely, we will show
that the classifier $\eta_n$, which is introduced in the next
section on the basis of a Transformer encoder, 
satisfies
\begin{equation}
  \label{se1eq3}
  \PROB\{ \eta_n(X) \neq Y \}
  -
  \min_{f:\R^{d \cdot l} \rightarrow \{0,1\}}  \PROB\{ f(X) \neq Y \}
  \leq
  c_1 \cdot (\log n)^{3} \cdot
  \max_{(p,K) \in \P} n^{- \frac{p}{2p+K}}
  \end{equation}
provided the aposteriori probability satisfies
a hierarchical
  composition model with some finite level and smoothness and order constraint
  $\P$. Since the above rate of convergence does not depend on the dimension
  $d \cdot l$ of $X$, our classifier is able to circumvent the curse of dimensionality
  if the aposteriori probability satisfies a suitable 
hierarchical
composition model.

In order to show (\ref{se1eq3}) we derive new approximation properties
and generalization bounds for Transformer decoders. The main idea
here is to show that the combination of attention units with
piecewise feedforward neural networks enables us to reconstruct
piecewise polynomials by Transformer decoders, and to generalize
a bound on the VC dimension of deep neural networks from
Bartlett et al. (2019) such that it is applicable to Transformer
encoders.
The Transformer classifiers
analyzed theoretically in this paper
are different from the one used in practice. We illustrate this
by describing classification problems in natural language processing
and the methods nowadays used in practice to solve these problems.
This gives us useful hints in which way the theoretical results
in this paper should be generalized in future work.

\subsection{Notation}
Throughout this paper, the following notation is used:
The sets of natural numbers, natural numbers including $0$,
integers
and real numbers
are denoted by $\N$, $\N_0$, $\Z$ and $\R$, respectively.
For $z \in \R$, we denote
the smallest integer greater than or equal to $z$ by
$\lceil z \rceil$.
Furthermore we set $z_+=\max\{z,0\}$.
Let $D \subseteq \R^d$ and let $f:\R^d \rightarrow \R$ be a real-valued
function defined on $\R^d$.
We write $x = \arg \min_{z \in D} f(z)$ if
$\min_{z \in \D} f(z)$ exists and if
$x$ satisfies
$x \in D$ and $f(x) = \min_{z \in \D} f(z)$.
For $f:\R^d \rightarrow \R$
\[
\|f\|_\infty = \sup_{x \in \R^d} |f(x)|
\]
is its supremum norm, and the supremum norm of $f$
on a set $A \subseteq \R^d$ is denoted by
\[
\|f\|_{A,\infty} = \sup_{x \in A} |f(x)|.
\]
Furthermore we define the norm $\| \cdot \|_{C^q(A)}$ of the smooth function space $C^q(A)$ by
\begin{align*}
\|f\|_{C^q(A)} :=\max\left\{\|\partial^{\bj}f\|_{\infty, A}: \|\bj\|_1 \leq q, \bj \in \N^d\right\}
\end{align*}
for any $f \in C^q(A)$, where
\[
\partial^{\bj}f = \frac{\partial^{j_1+\dots+j_d} f}{
\partial x_1^{j_1} \dots \partial x_d^{j_d}
}
\quad
\mbox{for }
\bj=(j_1, \dots,j_d)^T \in \N^d.
\]

Let $\F$ be a set of functions $f:\Rd \rightarrow \R$,
let $x_1, \dots, x_n \in \Rd$ and set $x_1^n=(x_1,\dots,x_n)$.
A finite collection $f_1, \dots, f_N:\Rd \rightarrow \R$
  is called an $\varepsilon$-- cover of $\F$ on $x_1^n$
  if for any $f \in \F$ there exists  $i \in \{1, \dots, N\}$
  such that
  \[
\frac{1}{n} \sum_{k=1}^n |f(x_k)-f_i(x_k)| < \varepsilon.
  \]
  The $\varepsilon$--covering number of $\F$ on $x_1^n$
  is the  size $N$ of the smallest $\varepsilon$--cover
  of $\F$ on $x_1^n$ and is denoted by $\Nu_1(\varepsilon,\F,x_1^n)$.

For $z \in \R$ and $\beta>0$ we define
$T_\beta z = \max\{-\beta, \min\{\beta,z\}\}$. 
For $i,j \in \N_0$ we let $\delta_{i,j}$ be the Kronecker delta, i.e.,
we set
\[
\delta_{i,j}=
\begin{cases}
  1, & \mbox{if } i=j,\\
  0, & \mbox{else}.
  \end{cases}
\]
If $W$ is a matrix and $b$ is a vector then we denote the number of
nonzero components in $W$ and $b$ by $\|W\|_0$ and $\|b\|_0$,
respectively.

\subsection{Outline}
In Section \ref{se2} our classifier
based on a Transformer encoder is defined.
The main result is presented in Section \ref{se3} and proven
in Section \ref{se5}. Section \ref{se4} describes the
application of Transformer classifiers in natural language
processing.

\section{Definition of a classifier based on a Transformer encoder}
\label{se2}

The Transformer encoder which we introduce in this
section becomes as input a sequence
\begin{equation}
  \label{se2eq1}
x=(x_1, \dots, x_l) \in \R^{l \cdot d}
\end{equation}
of length $l$ consisting of components $x_j \in \R^{d}$
$(j=1, \dots, l)$, where
$l$ and $d$ are natural numbers. As a first step
it produces from this sequence a new representation
\begin{equation}
  \label{se2eq3}
z_0=(z_{0,1}, \dots, z_{0,l}) \in \R^{l \cdot d_{model}}
\end{equation}
for some $d_{model} \in \N$. This new representation
is defined as follows:
We choose
$h,I \in \N$ and
set
\begin{equation}
  \label{se5eq3}
  d_{model}
  =
  h \cdot I \cdot (d+l+4).
\end{equation}
Here we repeat a coding of the input which includes the
original data, a coding of the position and additional
auxiliary values used for
later computation of function values
 $h \cdot I$ times.
More precisely we set for all $k \in \{1, \dots, h \cdot I\}$
\[
z_{0,j}^{((k-1) \cdot (d+l+4) +s)}
=
\begin{cases}
  x_j^{(s)} & \mbox{if } s \in \{1, \dots, d\} \\
  1 & \mbox{if } s=d+1 \\
  \delta_{s-d-1,j} & \mbox{if } s \in \{d+2, \dots, d+1+l\}\\
  1 & \mbox{if } s=d+l+3 \\
  0 & \mbox{if } s \in \{d+l+2, d+l+4 \}\\
  \end{cases}
\]

After that it computes successively
representations
\begin{equation}
  \label{se2eq2}
z_r=(z_{r,1}, \dots, z_{r,l})  \in \R^{l \cdot d_{model}}
\end{equation}
of the input for $r=1, \dots, N$, and uses $z_N$ in order
to predict a value $y \in \{0,1\}$. Here $N$ is the number
of pairs of attention and pointwise feedforward layers
of our Transformer encoder.

Given $z_{r-1}$ for some $r \in \{1, \dots, N\}$ we
compute $z_r$ by applying first a multi-head attention
and by applying second pointwise a feedforward neural network
with one hidden layer. Both times we will use an additional
residual connection.

The computation of the multi-head attention depends on matrices
\begin{equation}
  \label{se2eq4}
  W_{Q,r,s}, W_{K,r,s} \in \R^{d_k \times d_{model}}
    \quad \mbox{and} \quad
    W_{V,r,s} \in \R^{d_v \times d_{model}}
    \quad (s=1, \dots, h),
\end{equation}
where $h \in \N$ is the number of attentions which we compute
in parallel, where $d_k \in \N$ is the dimension of the queries
and the keys, and where $d_v=d_{model} /h= I \cdot (d+l+4)$
is the dimension of the values.
 We use these matrices to compute
for each component $z_{r-1,i}$ of $z_{r-1}$ corresponding queries
\begin{equation}
  \label{se2eq5}
  q_{r-1,s,i} = W_{Q,r,s} \cdot z_{r-1,i},
\end{equation}
keys
\begin{equation}
  \label{se2eq6}
  k_{r-1,s,i} = W_{K,r,s} \cdot z_{r-1,i}
\end{equation}
and values
\begin{equation}
  \label{se2eq7}
  v_{r-1,s,i} = W_{V,r,s} \cdot z_{r-1,i}
\end{equation}
$(s \in \{1, \dots, h\}, i \in \{1, \dots, l\})$.
Then the so-called attention between the component
$i$ of $z_{r-1}$ and the component $j$ of $z_{r-1}$ is defined
as the scalar product
\begin{equation}
  \label{se2eq8}
<q_{r-1,s,i}, k_{r-1,s,j}>
\end{equation}
and the index $j$ for which the maximal value occurs, i.e.,
\begin{equation}
  \label{se2eq9}
  \hat{j}_{r-1,s,i} = \argmax_{j \in \{1, \dots, l\}}
  <q_{r-1,s,i}, k_{r-1,s,j}>,
\end{equation}
is determined. The value corresponding to  this index is multiplied with the
maximal attention in (\ref{se2eq8}) in order to define
\begin{eqnarray}
  \label{se2eq10}
  \hspace*{-0.5cm}
  \bar{y}_{r,s,i}&=&v_{r-1,s,\hat{j}_{r-1,s,i}} \cdot \max_{j \in \{1, \dots, l\}}
  <q_{r-1,s,i}, k_{r-1,s,j}> \nonumber \\
  \hspace*{-0.5cm}
  &=& v_{r-1,s,\hat{j}_{r-1,s,i}} \cdot  <q_{r-1,s,i}, k_{r-1,s,\hat{j}_{r-1,s,i}}>
  \quad (s \in \{1, \dots, h\}, i \in \{1, \dots, l\}).
\end{eqnarray}
Using a residual connection 
we  compute the output of the multi-head attention
by
\begin{equation}
  \label{se2eq10}
  y_{r}=z_{r-1}+(\bar{y}_{r,1}, \dots, \bar{y}_{r,l})
\end{equation}
where
\[
\bar{y}_{r,s}
= (\bar{y}_{r,s,1}, \dots, \bar{y}_{r,s,h})
\in \R^{d_v \cdot h} = \R^{d_{model}}
\quad (s \in \{1, \dots, l\}).
\]
Here $y_r \in \R^{d_{model} \cdot l}$ has the same dimension as $z_{r-1}$.

The output of the pointwise feedforward neural network depends
on parameters
\begin{equation}
  \label{se2eq11}
  W_{r,1} \in \R^{d_{ff} \times d_{model}}, b_{r,1} \in \R^{d_{ff}},
  W_{r,2} \in \R^{d_{model} \times d_{ff}}, b_{r,2} \in \R^{d_{model}},
\end{equation}
which describe the weights in a feedforward neural network with
one hidden layer and $d_{ff} \in \N$ hidden neurons. This feedward
neural network is applied to each component of (\ref{se2eq10})
(which is analoguous to a convolutionary neural network) and computes
\begin{equation}
  \label{se2eq12}
  z_{r,s}=y_{r,s}+W_{r,2} \cdot \sigma \left(
  W_{r,1} \cdot y_{r,s} + b_{r,1}
  \right) + b_{r,2} \quad (s \in \{1, \dots, l\}),
\end{equation}
where we use again a residual connection.
Here
\[
\sigma(x) = \max\{x,0\}
\]
is the ReLU activation function, which is applied
to a vector by applying it
in each component
of the vector separately.

Given the output $z_N$ of the sequence of $N$ multi-head
attention and pointwise feedforward layers, our final
classifier is computed by
\[
\hat{y}
=
\begin{cases}
  1 & \mbox{if } z_N \cdot w + b \geq 1/2, \\
  0 & \mbox{else},
  \end{cases}
\]
where $w \in \R^{d_{model} \cdot l}$ and $b \in \R$ are parameters
of our neural network. I.e., we compute a linear transformation
of $z_N$ and use a plug-in classification rule corresponding to the function
\[
(x_1, \dots, x_l) \mapsto  z_N \cdot w + b.
\]
This function depends on a parameter vector
                        \[
                         \bar{\theta} = \left(
  (W_{Q,r,s}, W_{K,r,s}, W_{V,r,s})_{r \in \{1, \dots, N\}, s \in \{1, \dots, h\}},
    (W_{r,1},b_{r,1},W_{r,2},b_{r,2})_{r \in \{1, \dots, N\}}, w, b
                         \right).
                         \]
                         We denote this function by $f_{\bar{\theta}}$, i.e.,
                         $f_{\bar{\theta}}: \R^{d \cdot l} \rightarrow \R$ is the function
                         \[
f_{\bar{\theta}}(x_1, \dots, x_l)=  z_N \cdot w + b,
\]
where $z_N$ is the value computed as described above on the basis
of the matrices and vectors contained in $\bar{\theta}$.

In order to learn this function from observed data
$(X_1,Y_1)$, \dots, $(X_n,Y_n)$, we use the principle of least
squares to fit such a function to the observed data under a
sparsity constraint. To do this, we denote the number of
nonzero parameters of our transformer network by
\begin{eqnarray*}
\|\bar{\theta}\|_0
&=&
\sum_{r=1}^N  \sum_{s=1}^l \left(
\|W_{Q,r,s}\|_0
+
\|W_{K,r,s}\|_0
+
\|W_{V,r,s}\|_0
\right)
\\
&&
+
\sum_{r=1}^N
\left(
\|W_{r,1}\|_0
+
\|b_{r,1}\|_0
+
\|W_{r,2}\|_0
+
\|b_{r,2}\|_0
\right)
+
\|w\|_0
+
\|b\|_0,
\end{eqnarray*}
where $\| \cdot \|_0$ is the number of nonzero entries
in a matrix or in a vector. We then choose a sparsity
index $L_n \in \N$ and define
\begin{equation}
  \label{se2eq13}
  m_n(\cdot)
  =
  \arg \min_{f \in \{f_{\bar{\theta}} \, : \,
    \|\bar{\theta}\|_0 \leq L_n \} }
  \frac{1}{n} \sum_{i=1}^n | Y_i - f(X_i)|^2
\end{equation}
and
\begin{equation}
  \label{se2eq14}
  \eta_n(x)
  =
\begin{cases}
  1 & \mbox{if } m_n(x) \geq 1/2, \\
  0 & \mbox{else}.
  \end{cases}
\end{equation}
In (\ref{se2eq13}) we assume for simplicity that the minima does
indeed exist. If this is not the case, our main result below also
holds for any estimate which minimizes the empirical $L_2$ risk
in (\ref{se2eq13}) only up to some additional term of order $1/\sqrt{n}$.

\section{Main result}
\label{se3}

Our main result is the following theorem, which gives an upper
bound on the difference between the misclassification probability
of the estimate introduced in Section \ref{se2} and the
optimal misclassification probability.

\begin{theorem}
  \label{th1}
  Let $A>0$.
  Let $(X,Y)$, $(X_1, Y_1)$, \dots, $(X_n,Y_n)$ be
  independent and identically distributed $[-A,A]^{d \cdot l} \times \{0,1\}$--valued
  random variables, and let \linebreak
  $\eta(x)= \PROB\{Y=1|X=x\}$ be the corresponding
  aposterio probability. Let $\P$ be a finite subset
  of $(0,\infty) \times \N$ and assume that $\eta$ satisfies a hierarchical
  composition model with some finite level and smoothness and order constraint
  $\P$. Set
  \[
  h= \max_{(p,K) \in \P} n^{K/(2p+K)},
  I=\log n, L_n= (\log n)^2 \cdot
  (\max_{(p,K) \in \P} p) \cdot (\max_{(p,K) \in \P} K)
  \cdot \max_{(p,K) \in \P} n^{K/(2  \cdot p+K)},
  \]
  choose
$N \in \N$ sufficiently large,
  $d_K \geq 2$, $d_{ff} \geq 2 \cdot h + 2$, and define the estimate
  $\eta_n$ by (\ref{se2eq13}) and (\ref{se2eq14}). Then we have
  for $n$ sufficiently large
  \[
  \PROB\{ \eta_n(X) \neq Y \}
  -
  \min_{f:\R^{d \cdot l} \rightarrow \{0,1\}}  \PROB\{ f(X) \neq Y \}
  \leq
  c_1 \cdot (\log n)^{3} \cdot
  \max_{(p,K) \in \P} n^{- \frac{p}{2p+K}}.
  \]
  \end{theorem}
\noindent
    {\bf Remark 1.}
It follows from the proof of Theorem \ref{th1} that it also holds
if we set $I=c_2$ for some sufficiently large $c_2 \in \N$, which
depends on the level of the hierarchical composition model and
on $\P$.

    \noindent
        {\bf Remark 2.} The structure of the Transformer encoder
        in Theorem \ref{th1}, i.e., the multi-head attention and 
        piecewise feedforward layers and the residual connection,
        is as it is proposed in Vaswani et al. (2017). But different
        from the model proposed there is our choice of the encoding
        of the input. It is an open problem whether a similar
        result as in Theorem \ref{th1} also holds if the coding
        of the input position is done as in Vaswani et al. (2017).

\noindent
    {\bf Remark 3.} The rate of convergence in Theorem \ref{th1}
    does not depend on the dimension $d \cdot l$ of the predictor
    variable, hence the Transformer encoder is able to circumvent
    the curse of dimensionality in case that the aposteriori
    probability satisfies a hierarchical composition model
    with suitable order and smoothness constraints.

\section{Classification problems in natural language processing}
\label{se4}

In this section, we describe some of the \secrev{traditional}{} text classification datasets \secrev{}{that have been traditionally used to evaluate neural models} along with the typical NLP models to solve these tasks, including the Transformer based classifier. Then, we discuss the evaluation and the performance of the models on these datasets, and point out the differences between the
Transformer classifiers theoretically analyzed in this paper
and the ones used in practice. This gives us hints for useful
generalizations of the theoretical results in this paper in future
work.

%\subsection{Classification Task}

\subsection{Datasets}

The datasets have been constructed in a previous study by
Zhang, Zhao and LeCun (2015)
to foster the progress in text classification tasks and to empirically show that convolutional neural networks can provide competitive results on such tasks. 

\paragraph{YELP} contains reviews written by the customers about their experience in a location such as restaurants, bars and doctors. The dataset defines the classification task as predicting the number of stars ranging from 1 to 5 given the review text. It is a balanced dataset containing \rev{}{equal amount per \rev{star}{label} of training and testing samples} (130K training, 10K testing). This is known as Yelp F. (Full) where the task
is to predict the full stars. There is also a polarity classification version of this dataset, where 1 and 2 star reviews are merged under the Negative label, and 4 and 5 stars are labeled as Positive. This version is referred to as Yelp P. (Polarity). 

\paragraph{DBPedia} is a crowdsourcing project that maps Wikipedia infoboxes to a shared ontology (e.g., hierarchical classes of objects) that contains around 320 classes. The dataset that is derived from DBPedia consists of titles and abstracts of Wikipedia articles and their corresponding classes in the ontology. The authors construct the dataset by randomly sampling 40K training and 5K test samples from 14 non-overlapping classes from the DBPedia 2014 release.

\subsection{Models}

We present three techniques to tackle the aforementioned multilabel classification problems: \secrev{an n-gram tfidf model}{a strong baseline model using tfidf}, \secrev{convolutional neural network}{the early state-of-the-art-model based on convolutional neural networks} and the transformer based model.  

% The bag-of-ngrams models are constructed by selecting the 500,000 most frequent n-grams (up to 5-grams) from the training subset for each dataset
\subsubsection{Baseline} As a baseline, we choose the traditonal n-gram tf-idf model as in
Zhang, Zhao and LeCun (2015). \textit{n-gram} is one of the most common NLP terms that has been traditionally used to model text sequences. Here, $n$ refers to \secrev{a number (typically between 1 to 5) that describe the number sequential units of text (typically words)}{the number of sequential text units (typically words)} extracted from running text such as newspaper, wikipedia or web articles. In very simple terms, given the sentence ``I love NLP'', 1-grams that can be extracted from this sentence are \{ 'I', 'love', 'NLP'\}, while 2-grams are \{ 'I love', 'love NLP' \}. This baseline starts with extracting the most frequent 500,000 n-grams from the training sets, where $n$ is between 1 and 5. \rev{Then a simple classifier model is built using n-gram TFIDF features.}{}
Then it defines
TF (term-frequency) by
\begin{equation}
  tf(t,d) = \frac{f_{t,d}}{\sum_{t'\in d} f_{t',d}}
\end{equation}
Here $f_{t,d}$ is a raw count of the term (in this case a specific $n-gram$) that occurs in document $d$, which is then normalized by the number of times other terms ($t'$) occur in the same document, i.e., by the number of different terms which occur in $d$. It can be interpreted as certain types of documents containing certain n-grams, such as a sports article having a higher $tf$ for the 2-gram ``football game''.
Furthermore it defines IDF~(Inverse Document Frequency)
by
\begin{equation}
  idf(t,D) = log \frac{N}{|\{d \in D : t \in d\}|}
\end{equation}
Here $N$ is the total number of documents in the training set $D$. The denominator refers to the number of documents that contains the term $t$. It can therefore be interpreted as the amount of information a term provides. For instance the 1-gram ``the'' would have a low $idf$ since it is a common word, while ``football'' would yield a higher $idf$ since it occurs only on a subset of documents.
Using these two concepts it introduces
{\it tfidf} as follows:
\begin{equation}
  tfidf(t,d,D) = tf(t,d) \times idf(t,D)
\end{equation}
The values from the above equation are then used as features to train a multinomial logistic regression model.

%The first model is based on convolutional neural networks, and 
\subsubsection{Character-level Convolutional Neural Networks}

\begin{figure*}
\centering
    \includegraphics[width=1.0\textwidth]{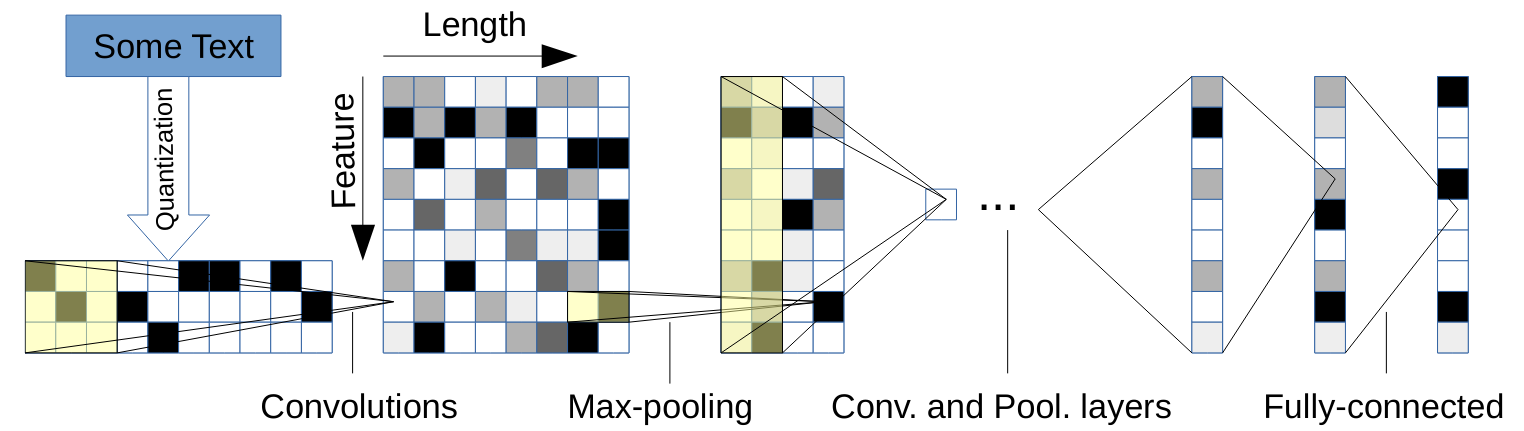}
    \caption{Character-level CNN architecture for text classification taken from
      Zhang, Zhao and LeCun (2015).}
  \label{fig:CNN}
\end{figure*}

The overall architecture of the character CNN is given in Fig.~\ref{fig:CNN}. First step is the \textbf{quantization} step, which is also referred to as the \textit{encoding} step, where characters are encoded as one-hot vectors, \rev{}{by marking the character to be represented with the value $1$ in a zero-vector (e.g., [0, 0, 1, 0, 0, .., 0])}. \rev{}{More formally, given} an alphabet of $k$ characters, each character is indexed and represented with a vector $c \in \R^k$, where the component corresponding to the number of the character in the alphabet
is $1$, and in all other components the value is $0$. In Zhang, Zhao and LeCun (2015) an alphabet of 70 characters is defined which is a combination of 26-letter English alphabet, digits, punctuations and the space character.

\paragraph{1D Convolution}: Given the sequence of characters $c_{1:n}=(c_1,...,c_n)$, a one dimensional convolution can simply be defined as sliding a window of size $k$ over the character sequence to apply a convoluton filter, a.k.a., kernel, to each window. More formally, given the filter $u=(u_1, \dots, u_k) \in \R ^ {d \times k}$, and the window of characters $x_i = (c_i, c_{i+1},...,c_{i+k}) \in \R ^ {d \times k}$, the convolution is defined as a dot-product between $x_i$ and $u$,
i.e., by
\[
u \cdot x_j
=
\sum_{j=1}^k u_j \cdot c_{i+j-1}.
\]
In practice, however, we use a series of $l$ filters $U_1, \dots, U_l \in \R ^ {d \times k}$ and a bias term $b=(b_1, \dots, b_l) \in \R^l$. The convolution is then followed by the nonlinear ReLU activation function $g(x) = max\{0, x\}$,
i.e., the result of the convolution is
\begin{equation}
  s_{i,j} = g(x_i \cdot U_j + b) \quad (i=1, \dots, n-k, j=1, \dots, l).
\end{equation}

\paragraph{Max Pooling}: This operation simply chooses the maximum of the values from all $s_i$ vectors, which yields the fixed $l$ dimensional vector $h=(h_1, \dots, h_l)^T$ where
\[
h_j = \max_{i=1, \dots, n-k} s_{i,j}.
\]
The sequence of convolution and max pooling operations are referred to as a \textbf{convolutional layer}. 

\rev{}{The convolution layer is then followed by a linear layer with weight matrix, $W$, and a bias term $b$. The output of the linear layer is then fed to a softmax function to calculate the label probabilities for each class, $l_i$, as follows:} 
\begin{equation}
  p(l_i) = softmax(W \cdot h + b).
\end{equation}

The model is trained with back propagation algorithm using stochastic gradient descent \rev{}{to minimize the cross-entropy loss, $L_{CE}$,}
\begin{equation}
  L_{CE} = -\sum_{i=1}^{n} {t_i \cdot log(p(l_i))}
\end{equation}
\rev{}{where $n$ refers to the number of classes, and $t$ is a one-hot vector of dimension $n$ where the true label of the sample is flagged with $1$ (e.g., if the test sample is from class 2, then $t=(t_1, \dots, t_n)$ will be $(0, 1, 0, .. , 0)$).}

The final model consists of 6 convolutional layers which are responsible \secrev{from}{of} extracting relevant textual features, which is then followed by 3 linear layers that capture the relation between the features and the classification labels. 

% More details on weight initialization etc...

\subsubsection{Transformer based Language Models}

One of the recent breakthroughs in natural language processing is introduction of pretrained large language models, a.k.a., LLMs, that use the Transformer architecture such as Bidirectional Encoder Representations from Transformers (BERT),
cf., Devlin et al. (2019). Such pretrained models also trigger a new paradigm for training \rev{downstream}{high-level} NLP tasks \rev{}{such as document classification, named entity recognition and question answering}. This new training regime, called \textbf{fine-tuning}, allows the NLP practitioners and researchers to benefit from unlabeled text via performing small weight modifications to LLMs by adding a task-specific output layer and performing a minimal supervised re-training on the task data. 
\gozde{Since BERT has been the predominantly used LLM, we experiment with BERT in this paper. Let us explain the architecture, training objective and the fine-tuning process in more details.} 

% The advantage of these approaches is that few parameters need to be learned from scratch

% what is BERT

% ELMO - separate training of bidirectional objective and then concatenate the two representations
\gozde{
  \paragraph{Architecture} The architecture of the BERT model is identical to the multi-layer bidirectional Transformer encoder described in Vaswani et al. (2017) and released in \textit{tensor2tensor} library\footnote{https://github.com/tensorflow/tensor2tensor}.} This means that it corresponds exactly to
the Transformer encoder defined in Section \ref{se2}.

\paragraph{Training Objective} Traditional language modeling objective has been defined as \textit{predicting the next token given the previous context}. More formally, the objective has been to maximize the log probability $\log P(w_{i+1}|w_i, w_{i-1},...,w_1)$ where $w_{i}$ denotes the word at index $i$. Later, researchers introduced a secondary objective which is to \textit{predict the previous token given the future context}, i.e., to maximize $\log P(w_{i-1}|w_i, w_{i+1},...,w_n)$ where the maximum number of words is $n$. Traditionally, two different representations are learned using these two objectives which are then concatenated and used as the final token representation. 

% simply mask some percentage of the input tokens at random
The calculation of this objective has been trivial when sequential neural models such as recurrent neural networks have been used. However, since Transformer architecture is not designed to be autoregressive, such bidirectional training objective has not been that intuitive. To address this gap, BERT uses a combination of different language modeling objectives. \gozde{The first objective is called a \textbf{Masked Language Model (MLM)}. It is defined as predicting the tokens that are masked at random. To learn representations that can better model the relation between sentences, BERT defines the second self-explanatory objective: \textbf{Next Sentence Prediction (NSP)}.} 

In BERT, words are represented as a sequence of subwords using the WordPiece algorithm which is inspired from an old byte compression idea. The authors introduce a special token named \textsc{[CLS]} that is placed at the beginning of each sequence. The final hidden state of this token \gozde{from the multi-layer bidirectional Transformer encoder}, i.e., the hidden state of the pointwise
feedforward neural network in
(\ref{se2eq12})
corresponding to this token,
is considered as the sequence representation and used for the classification tasks. Since BERT aims to provide a unified architecture for NLP downstream tasks, including tasks that take a \textit{pair} of inputs (e.g., question and answer), it introduces another special token called \textsc{[SEP]} to separate two sentences from each other. The input representation of BERT is then the sum of three embeddings: the token (or the WordPiece) embedding, \textbf{$E_{tok}$}, the segment embedding \textbf{$E_{A|B}$} that identifies whether the sentence is the prior (segment A) or the next (segment B) sentence, and a special \textbf{Positional Embedding (PE)} using sinus and cosine wave frequencies \rev{}{defined as follows:}
\begin{equation}
PE_{pos, 2i} = sin \frac{pos}{10000^{\frac{2i}{d}}},  \qquad PE_{pos, 2i+1} = cos \frac{pos}{10000^{\frac{2i}{d}}}, 
\end{equation}
\rev{}{where $d$ is the fixed dimension size of the input representation, $i$ refers to the index \gozde{in the $d$ dimensional Positional Embedding (PE) and $pos$ refers to the position of the token within the input. For instance the $pos$ value of the token ``I'' in the input ``I love NLP'' would be 0, and ``love'' would be 1.}} 

% the final representation: token + position + segment

\paragraph{Pretraining} The model is pretrained with the two objectives, \textbf{Masked Language Model (MLM)} and \textbf{Next Sentence Prediction (NSP)}, that are mentioned above. \gozde{To elaborate,} for MLM, 15\% of the WordPiece tokens are masked randomly. Later the hidden vectors for the masked tokens are passed onto a softmax layer that calculates the probability distribution over the full vocabulary. NSP is formulated as a binary classification problem, where \textsc{[CLS]} token is classified either as \textsc{isNext} or \textsc{notNext}. It can be interpreted as given two sentences, the second sentence is either the \textit{next} sentence following the previous one, or not. The training loss is then the sum of the mean MLM and NSP likelihood. The model is pretrained on a corpus of a large collection (in total 3,3 billion words) of books and wikipedia articles. It should be noted that, in practice, sequence pairs do not directly refer to linguistically well-formed sentence structures. \rev{}{For instance, a sequence might contain an unfinished sentence since it contains a fixed number of tokens and does not take sentence boundaries into account.} They are typically longer spans of text which may contain multiple sentences. For training proper segment embeddings, training data is constructed in a way that, 50\% of the time the second sentence follows the first sentence (i.e., has the label isNext) and 50\% it is a randomly chosen sentence (i.e., has the label notNext). %Since training longer sequences are extremely expensive, the model is first trained with sequences of length 128, and then   

\begin{figure*}
\centering
    \includegraphics[width=1.0\textwidth]{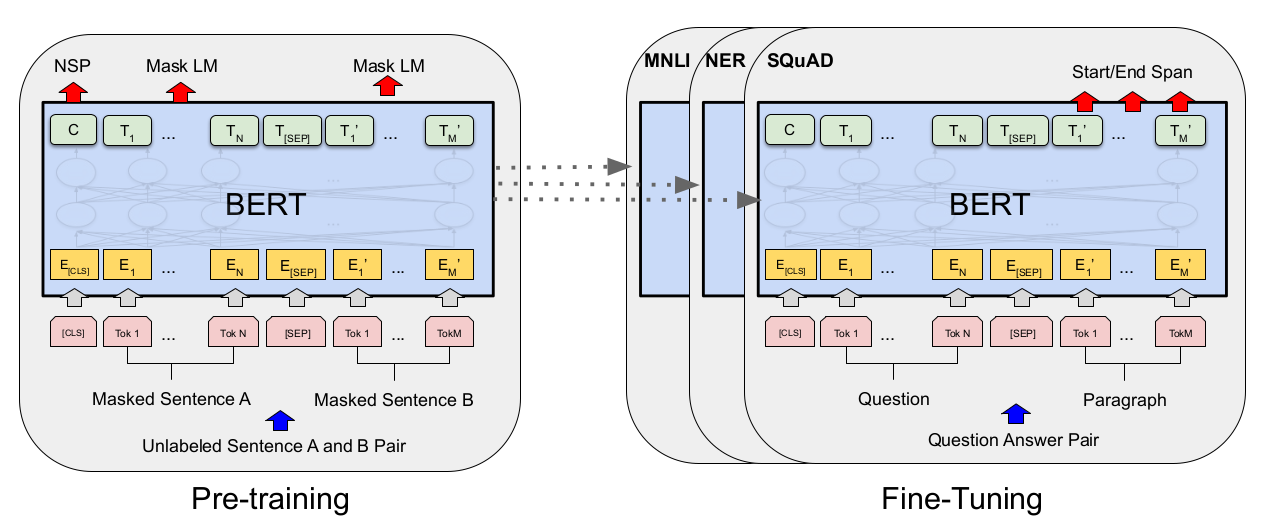}
    \caption{Illustration of pretraining and finetuning of BERT taken from
      Devlin et al. (2019). $E$ refers to embedding, \textsc{[CLS]} and \textsc{[SEP]} denote special tokens designed to perform classification and to separate two sentences, accordingly. NSP denotes ``next sentence prediction''. MNLI~(Multi-Genre Natural Language Inference), NER~(Named Entity Recognition) and SQuAD~(Stanford Question Answering Dataset) are popular NLP tasks for on which the BERT model is finetuned.}
  \label{fig:BERT}
\end{figure*}

\paragraph{Fine-tuning} Fine-tuning process, in general, aims to introduce a minimal amount of task-specific parameters which
are used to define a task-specific output layer. This output layer depends on the \rev{downstream}{NLP} task. For instance, for token-level tasks such as \rev{NER}{Named Entity Recognition (NER)} (i.e., labeling individual tokens as entity names such as organization, company and city name), token representations are fed into the output layer. In case of classification tasks \rev{such as sentiment analysis}{such classifying restaurant reviews as in this paper},
the final hidden state of the \textsc{[CLS]} token is fed into the output layer. Finally,
either all parameters or only the task-specific parameters are updated end-to-end using the task-specific objective \rev{}{(e.g., binary cross-entropy loss for binary classification, categorical cross-entropy loss for multiple classes)} on the \rev{}{task-specific} labeled data \rev{}{(e.g., the classification datasets described in Sec.~4.1)}. An illustration of the pretraining and the fine-tuning processes is given in Fig.~\ref{fig:BERT}. 

% BERT alleviates the previously mentioned unidirectionality constraint by using a “masked language model” (MLM) pre-training objective, inspired by the Cloze task (Taylor, 1953).    

% The masked language model randomly masks some of the tokens from the input, and the objective is to predict the original vocabulary id of the masked word based only on its context.

% the MLM objective enables the representation to fuse the left and the right context, which allows us to pretrain a deep bidirectional Transformer.

% In addition to the masked language model, we also use a “next sentence prediction” task that jointly pretrains text-pair representations.

% The fine-tuning approach, such as the Generative Pre-trained Transformer (OpenAI GPT) (Radford et al., 2018), introduces minimal task-specific parameters, and is trained on the downstream tasks by simply fine-tuning all pretrained parameters.

\subsection{Evaluation and Results}
All models are evaluated on the \rev{}{standard} test splits \rev{}{as provided by the data created in Zhang, Zhao and LeCun (2015)} using the empirical
misclassification error on the testing data as the evaluation measure. We show the scores of the baseline as reported by Zhang, Zhao and LeCun (2015). 

\secrev{For BERT, there exists many additional finetuning strategies. Sun et al.
  (2019)  study such engineering choices for text classification tasks and report three efficient finetuning strategies: a) Within task pretraining: BERT is further pre-trained on the training data of a target task, b) In-Domain pretraining: further pretraining on a combination of training instances of related tasks which are assumed to be from the same domain, c) Cross-domain pretraining: further pretraining on a mixture of domains. We show the results of these models as reported in Sun et al.
  (2019) in Table~\ref{tab:results}.}
       {We finetune the BERT-base model that consists of 12 Transformer encoder blocks with hidden dimensionality of 768, and 12 self-attention heads. BERT outputs a representation of an input sequence including a representation for the special \textsc{[CLS]} token. The maximum number of tokens that can be processed by BERT is 512. However, since processing 512 tokens is quite memory demanding for a single GPU, we set the threshold to 128 tokens. We take the final hidden state $h$ of the \textsc{[CLS]} token and introduce a task specific weight matrix $W$. The output is then fed to a softmax classifier to calculate the probability
of class $l$ by        
\begin{equation}
  p(l|h) = softmax(Wh).
\end{equation}
BERT has been shown to capture lower-level linguistic information like sentence structure in lower layers, while higher-level layers capture higher-level information such as the sentence semantics. Therefore, instead of fine-tuning all layers, we fine-tune only the last two layers. This strategy does not only reduce the computational workload, but also has been shown to provide competitive scores, cf., Sun et al.
  (2019). Finally, we maximize the log-probability of the correct label by updating the parameters in $W$ and the last two layers. We fine-tune the uncased model (all lower case) for 3 epochs with the batch size as 32 and the learning rate as 0.00002. We train \rev{}{both the character-level CNN and the fine-tuned BERT models} using five different random seeds and report the averaged scores in Table~\ref{tab:results}.} 

\secrev{As can be seen, the best results are consistently achieved by a BERT variant. The error rates reported from a BERT variant are \textit{marginally} better than the baseline and the CNN  model, having around 8\% points gap on some datasets such as Yelp F. It can be seen that, even a basic fine-tuning strategy shown with BERT-FiT overperforms other models by a large margin on all datasets.}{Our results show that, even a simple fine-tuning strategy on the large pretrained Transformer-based language model BERT provides substantial improvements over the baseline and the character-level CNNs.}

\begin{table*}[!htp]\centering
    %\scriptsize
    \scalebox{0.9}{
    \begin{tabular}{lccc}
    \toprule
    Model & Yelp P. & Yelp F. & DBP \\
    \midrule
    \textbf{n-gram TFIDF~(Zhang, Zhao and LeCun (2015))} & 4.56 & 45.20 & 1.31 \\
    %\textbf{Char CNN~(Zhang et. al., 2015 \cite{CNN15})} & 4.88 & 37.95 & 1.55\\
    \textbf{Char CNN} & 4.96 & 39.25 & 1.78\\
    \textbf{BERT Fine-Tuning} & 2.92 & 32.94 & 0.75 \\
    %\textbf{BERT-ITPT-FiT (Sun et. al., 2019)} & 1.92 & 29.42 & 0.68 \\
    %\textbf{BERT-IDPT-FiT (Sun et. al., 2019)} & \textbf{1.87} & 29.25 & \textbf{0.65} \\
    %\textbf{BERT-CDPT-FiT (Sun et. al., 2019)} & 1.97 & \textbf{29.20} & 0.67 \\
    \bottomrule
    \end{tabular}
    }
    \caption{Comparison of test error rates of the baseline, character-level CNN and the finetuned BERT model.}
    \label{tab:results}
  \end{table*}

% put a result table, add a few lines of discussion

\subsection{Comparison of theory and practice}
If we compare our classifier introduced in Section \ref{se2}
with the ones used in practice, we see that the basic
structure of the classifiers, i.e., the use of a sequence of pairs
of mulit-head attention and pointwise feedforward layers, is the
same. Different are the coding of the input (here we used a separate
part of the vector in order to encode the position instead of adding
some kind of encoded wave to the vector) and the learning method.
In particular, our theory totally ignores the idea of using a pretrained
model as used in BERT. Instead we focused only on the approximation and
generalization properties of Transformer networks in order to simplify
the theory.

This indicates what the main question for
future research in this area is: Can one also come up with a theory
for Transformer classifiers which uses firstly pretraining and
secondly learning of the classifiers by gradient descent?

\section{Proofs}
\label{se5}

\subsection{Proof of Theorem \ref{th1}}
Set $\beta=1$. Since $T_\beta z \geq 1/2$ if and only if $z \geq 1/2$
for any $z \in \R$ we have
\[
  \eta_n(x)
  =
\begin{cases}
  1 & \mbox{if } T_\beta m_n(x) \geq 1/2, \\
  0 & \mbox{else}.
  \end{cases}
\]
By Theorem 1.1 in Gy\"orfi et al. (2002) we know
\[
\PROB\{\eta_n(X) \neq Y|\D_n\}
-
  \min_{f:\R^{d \cdot l} \rightarrow \{0,1\}}  \PROB\{ f(X) \neq Y \}
\leq
2 \cdot
\sqrt{
\int |T_\beta m_n(x)-\eta(x)|^2  \PROB_X(dx)
  },
\]
hence it suffices to show
\begin{equation}
  \label{pth1eq1}
\EXP
\int |T_\beta m_n(x)-\eta(x)|^2  \PROB_X(dx)
\leq
  c_3 \cdot (\log n)^6 \cdot
  \max_{(p,K) \in \P} n^{- \frac{2p}{2p+K}}.
\end{equation}
By standard error bounds from empirical process theory
(cf., e.g.,  Lemma 1 in Bauer and Kohler (2019))
we know
    \begin{eqnarray*}
        &&
        \EXP \int |T_{\beta} m_n(x)-\eta(x)|^2  \PROB_X(dx)
        \\
        &&
        \leq
        \frac{c_{4} \cdot \left(\log n )^2 \cdot
          \sup_{x_1^n }
          \left(
          \log
          \left(
          \mathcal{N}_1
          \left(
          \frac{1}{n\cdot \beta}, 
\{ T_{\beta} f_{\bar{\theta}} \, : \,
    \|\bar{\theta}\|_0 \leq L_n \}
    \right), x_1^n
    \right)
    \right)
    +1
    \right)}{n} \\
&&\quad + 2 \cdot \inf_{f \in \{f_{\bar{\theta}} \, : \,
    \|\bar{\theta}\|_0 \leq L_n \}} \int |f(x)-\eta(x)|^2 {\PROB}_X (dx).
        \end{eqnarray*}
    We will show in Lemma \ref{le6} below that we have
    \begin{eqnarray*}
    \sup_{x_1^n } \log
    \left(
    \mathcal{N}_1
    \left(\frac{1}{n\cdot \beta}, 
\{ T_{\beta} f_{\bar{\theta}} \, : \,
    \|\bar{\theta}\|_0 \leq L_n \}, x_1^n
    \right)
    \right)
    &\leq&  c_5 \cdot L_n \cdot (\log n)^2
    \\
    &=&
    c_6
    \cdot (\log n)^4 \cdot \max_{(p,K) \in \P} n^{K/(2 \cdot p+K)}.
    \end{eqnarray*}
    And by Theorem \ref{th2} below we get
    \[
\inf_{f \in \{f_{\bar{\theta}} \, : \,
  \|\bar{\theta}\|_0 \leq L_n \}} \int |f(x)-\eta(x)|^2 {\PROB}_X (dx)
\leq c_7 \cdot   \max_{(p,K) \in \P} h^{-2 \cdot p/K}.
\]
Choose $(\bar{p},\bar{K}) \in \P$ with
\[
\frac{\bar{p}}{\bar{K}}
=
\min_{(p,K) \in \P} \frac{p}{K}.
\]
Then the above results imply
\begin{eqnarray*}
\EXP
\int |T_\beta m_n(x)-\eta(x)|^2  \PROB_X(dx)
&
\leq &
\frac{    c_6
    \cdot (\log n)^6 \cdot n^{\bar{K}/(2 \cdot \bar{p}+\bar{K})}
}{n}
+
c_7 \cdot    h^{-2 \cdot \bar{p}/\bar{K}}
\\
&
\leq &
c_8
\cdot (\log n)^6 \cdot n^{-2 \cdot \bar{p}/(2 \cdot \bar{p}+\bar{K})}
\\
&=&
c_8
\cdot (\log n)^6 \cdot
  \max_{(p,K) \in \P} n^{- \frac{2p}{2p+K}}.
  \end{eqnarray*}
    \hfill $\Box$

\subsection{Approximation of piecewise polynomials by transformer networks
  with single-head attention}
In this subsection we present results concerning the approximation
of piecewise polynomials by Transformer encoders
with single-head attention. In the next subsection we will
generalize these results to Transformer encoders
with multi-head attention.

In case of a Transformer encoder with a single-head attention
we represent the input sequence (\ref{se2eq1}) by
\[
z_0=(z_{0,1}, \dots, z_{0,l}) \in \R^{l \cdot d_{model}},
\]
where
\begin{equation}
  \label{se5eq1}
  d_{model}
  =
  d+l+4.
\end{equation}
Here we use $d$ components to represent the original data,
we use $l$ components to represent the position in the original
sequence, we use 1 component to have available the constant 1
in our computations, and we use 3 components to be able to
save auxiliary results during our computation. More precisely we set
\[
z_{0,j}^{(s)}
=
\begin{cases}
  x_j^{(s)} & \mbox{if } s \in \{1, \dots, d\} \\
  1 & \mbox{if } s=d+1 \\
  \delta_{s-d-1,j} & \mbox{if } s \in \{d+2, \dots, d+1+l\}\\
  1 & \mbox{if } s=d+l+3 \\
  0 & \mbox{if } s \in \{d+l+2, d+l+4 \}\\
  \end{cases}
\]

We start with a first lemma that shows that a single-head attention
unit can be used to compute linear polynomials in one variable.

\begin{lemma}
  \label{le1}
  Let $x_j \in \R^d$ and $b_j \in \R$ $(j=1, \dots, l)$.
  Let $z_0=(z_{0,1}, \dots, z_{0,l}) \in \Rd$ be given by
\[
z_{0,j}^{(s)}
=
\begin{cases}
  x_j^{(s)} & \mbox{if } s \in \{1, \dots, d\} \\
  1 & \mbox{if } s=d+1 \\
  \delta_{s-d-1,j} & \mbox{if } s \in \{d+2, \dots, d+1+l\}\\
  b_j & \mbox{if } s=d+l+3 \\
  0 & \mbox{if } s \in \{d+l+2, d+l+4 \}\\
  \end{cases}
\]
Let $j \in \{1, \dots, l\}$, $k \in \{1, \dots, d\}$ and $u \in \R$
be arbitrary. Let
\[
B > 2 \cdot \max_{s=1, \dots, d,j=1, \dots,l } |x_j^{(s)}|.
\]
Then there exists matrices
$  W_{Q}$, $W_{K} \in \R^{2 \times d_{model}}$ and
$W_{V} \in \R^{d_{model} \times d_{model}}$,
where each matrix contains at most 3 nonzero entries
and where all entries depend only on $u$ and $B$,
such that
  $q_{i} = W_{Q} \cdot z_{0,i}$,
  $k_{i} = W_{K} \cdot z_{0,i}$,
$v_{i} = W_{V} \cdot z_{0,i}$,
\[
  \hat{j}_{i} = \argmax_{t \in \{1, \dots, l\}}
  <q_{i}, k_{t}> \quad (i \in \{1, \dots, l\})
  \]
  and
  \[
  y_{i} = z_{0,i}+v_{\hat{j}_{i}} \cdot  <q_{i}, k_{\hat{j}_{i}}>
  \quad (i \in \{1, \dots, l\})
\]
result in
\[
y_{1}^{(s)}
=
\begin{cases}
  x_1^{(s)} & \mbox{if } s \in \{1, \dots, d\} \\
  1 & \mbox{if } s=d+1 \\
  \delta_{s-d-1,j} & \mbox{if } s \in \{d+2, \dots, d+1+l\}\\
  x_j^{(k)}-u  & \mbox{if } s= d+l+2\\
  b_1 & \mbox{if } s=d+l+3 \\
  0 & \mbox{if } s= d+l+4 \\
  \end{cases}
\]
and
\[
y_j=z_{0,j} \quad \mbox{for } j \in \{2, \dots, l \}.
\]
  \end{lemma}

\noindent
    {\bf Proof.} We start our proof by defining the matrices
    $W_Q$, $W_K$ and $W_V$. We set
    \[
    W_Q = \left(
    \begin{array}{ccccccc}
      0 & \dots & 0 & 1 & 0 & \dots & 0 \\
      0 & \dots & 0 & B & 0 & \dots & 0 \\
      \end{array}
    \right)
    \]
    where all colums are zero except column number $d+2$,
    \[
   W_K = \left(
    \begin{array}{ccccccccccccccc}
      0 & \dots & 0 & 1 & 0 & \dots & 0 & -u-B & 0 & \dots & 0 & 0 & 0 & \dots & 0 \\
      0 & \dots & 0 & 0 & 0 & \dots & 0 & 0 & 0 & \dots & 0 & 1 & 0 & \dots & 0 \\
      \end{array}
    \right)
    \]
    where all colums are zero except column number $k$, column number $d+1$
    and column number $d+1+j$, and
    \[
   W_V = \left(
    \begin{array}{ccccccc}
      0 & \dots & 0 & 0 & 0 & \dots & 0 \\
      \vdots & \dots & \vdots & \vdots & \vdots & \dots & \vdots\\
      0 & \dots & 0 & 0 & 0 & \dots & 0 \\
      0 & \dots & 0 & 1 & 0 & \dots & 0 \\
            0 & \dots & 0 & 0 & 0 & \dots & 0 \\
      \vdots & \dots & \vdots & \vdots & \vdots & \dots & \vdots\\
      0 & \dots & 0 & 0 & 0 & \dots & 0 \\
    \end{array}
    \right)
    \]
    where all rows and all colums are zero except row number $d+l+2$
    and column number $d+1$.
    For these matrices we get
    \[
    q_{i} = W_{Q} \cdot z_{0,i} =
    \left(
    \begin{array}{c}
      1 \cdot \delta_{i,1} \\
      B \cdot \delta_{i,1} \\
      \end{array}
    \right),
    \quad 
    k_{i} = W_{K} \cdot z_{0,i} =
    \left(
    \begin{array}{c}
      x_i^{(k)} -u -B  \\
      1 \cdot \delta_{i,j} \\
      \end{array}
    \right),
    \]
    \[
    v_{i} = W_{V} \cdot z_{0,i} = ( \delta_{1,d+l+2},\delta_{2,d+l+2}, \dots,
    \delta_{d+l+4,d+l+2})^T = e_{d+l+2},
    \]
    where $e_{d+l+2}$ is the $d+l+2$-th unit vector in $\R^{d+l+4}$.
    Consequently we have
    \[
    <q_r,k_s> \, = \, (x_s^{(k)}-u-B) \cdot \delta_{r,1} +B  \cdot \delta_{r,1}
    \cdot \delta_{s,j}, 
    \]
    \[
    <q_1,k_j> \, = \, (x_j^{(k)}-u-B)  +B \, > \, x_s^{(k)}-u-B \, = \, <q_1,k_s> 
    \]
    for all $s \in \{1, \dots, d\} \setminus \{j\}$
and
    \[
\hat{j}_1 = j,
\]
which implies
\[
v_{\hat{j}_{1}} \cdot  <q_{1}, k_{\hat{j}_{1}}>  = (x_j^{(k)}-u) \cdot
e_{d+l+2}
\]
and
\[
v_{\hat{j}_{i}} \cdot  <q_{i}, k_{\hat{j}_{i}}>  =
e_{d+l+2} \cdot  <0, k_{\hat{j}_{i}}>  =
0 \cdot
e_{d+l+2}
\]
for $i > 1$.
    \hfill $\Box$ \\

    \noindent
        {\bf Remark 4.}
        It follows from the proof of Lemma \ref{le1} that
        we can modify $W_K$ such that
        $y_1^{(d+l+2)}=1$ holds instead of $y_1^{(d+l+2)}=x_j^{(k)}-u$.\\

Our next  lemma  shows that a single-head attention
unit can be used to compute products.

\begin{lemma}
  \label{le2}
  Let $x_s \in \R^d$ and $a_s,b_s \in \R$ $(s=1, \dots, l)$.
  Let $z_0=(z_{0,1}, \dots, z_{0,l}) \in \Rd$ be given by
\[
z_{0,j}^{(s)}
=
\begin{cases}
  x_j^{(s)} & \mbox{if } s \in \{1, \dots, d\} \\
  1 & \mbox{if } s=d+1 \\
  \delta_{s-d-1,j} & \mbox{if } s \in \{d+2, \dots, d+1+l\}\\
  a_j & \mbox{if } s=d+l+2 \\
  b_j & \mbox{if } s=d+l+3 \\
  0 & \mbox{if } s=d+l+4 \\
  \end{cases}
\]
Let $j \in \{1, \dots, l\}$. Let $B >  2 \cdot \max_{r,s} |a_r \cdot b_s|$.
Then there exists matrices
$  W_{Q}$, $W_{K} \in \R^{2 \times d_{model}}$ and
    $W_{V} \in \R^{d_{model} \times d_{model}}$,
where each matrix contains at most 2 nonzero entries
and where all entries depend only on $B$,
such that
  $q_{i} = W_{Q} \cdot z_{0,i}$,
  $k_{i} = W_{K} \cdot z_{0,i}$,
$v_{i} = W_{V} \cdot z_{0,i}$,
\[
  \hat{j}_{i} = \argmax_{t \in \{1, \dots, l\}}
  <q_{i}, k_{t}> \quad (i \in \{1, \dots, l\})
  \]
  and
  \[
  y_{i} = z_{0,i}+v_{\hat{j}_{i}} \cdot  <q_{i}, k_{\hat{j}_{i}}>
  \quad (i \in \{1, \dots, l\})
\]
result in
\[
y_{1}^{(s)}
=
\begin{cases}
  x_1^{(s)} & \mbox{if } s \in \{1, \dots, d\} \\
  1 & \mbox{if } s=d+1 \\
  \delta_{s-d-1,j} & \mbox{if } s \in \{d+2, \dots, d+1+l\}\\
  a_1 & \mbox{if } s=d+l+2 \\
  b_1  & \mbox{if } s= d+l+3\\
  b_1 \cdot a_j + B & \mbox{if } s= d+l+4 \\
  \end{cases}
\]
and
\[
y_s^{(i)}=z_{0,s}^{(i)} \quad \mbox{for }
i \in \{1, \dots, d+l+3\},
s \in \{1, \dots, l\}. 
\]  
  \end{lemma}

\noindent
    {\bf Proof.} Define $W_V$ as in the proof of Lemma \ref{le1}
    such that
all rows and all colums are zero except row number $d+l+4$
    and column number $d+1$,    
    and
    set
    \[
    W_Q= \left(
    \begin{array}{ccccccccccc}
      0 & \dots & 0 & 0 &  0 &\dots & 0 & 1 &  0 & \dots & 0 \\
      0 & \dots & 0 & B &  0 &\dots & 0 & 0 &  0 & \dots & 0 \\
      \end{array}
    \right)
    \]
    where all colums are zero except column number $d+2$ and
    column number $d+l+3$ and
    \[
    W_K= \left(
    \begin{array}{ccccccccccc}
      0 & \dots & 0 & 0 &  0 &\dots & 0 & 1 &  0 & \dots & 0 \\
      0 & \dots & 0 & 1 &  0 &\dots & 0 & 0 &  0 & \dots & 0 \\
      \end{array}
    \right)
    \]
    where all colums are zero except column number $d+1+j$ and
    column number $d+l+2$.
    
    For these matrices we get
    \[
    q_{i} = W_{Q} \cdot z_{0,i} =
    \left(
    \begin{array}{c}
      b_i  \\
      B \cdot \delta_{i,1} \\
      \end{array}
    \right),
    \quad 
    k_{i} = W_{K} \cdot z_{0,i} =
    \left(
    \begin{array}{c}
      a_i   \\
      1 \cdot \delta_{i,j} \\
      \end{array}
    \right),
    \]
    \[
    v_{i} = W_{V} \cdot z_{0,i} = ( 0, \dots,0,1)^T = e_{d+l+4}.
    \]
    Consequently we have
    \[
    <q_r,k_s> = b_r \cdot a_s 
    +B   \cdot \delta_{r,1} \cdot \delta_{s,j}. 
\]
This implies
    \[
\hat{j}_1 = j
\]
and
\[
 <q_{i}, k_{\hat{j}_{i}}> = \max_{s} b_i \cdot a_s
 \]
 for $i>1$, from which we can conclude
\[
v_{\hat{j}_{1}} \cdot  <q_{1}, k_{\hat{j}_{1}}>  = (b_1 \cdot a_j + B) \cdot
e_{d+l+4}
\]
and
\[
v_{\hat{j}_{i}} \cdot  <q_{i}, k_{\hat{j}_{i}}>  =
(\max_{s} b_i \cdot a_s )
\cdot
e_{d+l+4}
\]
for $i > 1$.
    \hfill $\Box$

    Our next lemma defines a special pointwise feedforward neural network,
    which applies the function
    \[
x \mapsto \alpha \cdot (x-B)
    \]
    to  component number $d+l+4$ of each element in the decoding of the sequence
    of inputs and writes the result in component number $d+l+3$, and which
    sets the entries in components $d+l+2$ and $d+l+4$ to zero.

    \begin{lemma}
      \label{le3}
      Let $d_{model}=d+l+4$ and
      let $y=(y_1, \dots, y_l)$ for some $y_i \in \R^{d_{model}}$.
      Let $d_{ff} \geq 8$ and let $\alpha \in \R$.
      Then there exists matrices and vectors
      \[
      W_{1} \in \R^{d_{ff} \times d_{model}}, b_{1} \in \R^{d_{ff}},
      W_{2} \in \R^{d_{model} \times d_{ff}}, b_{2} \in \R^{d_{model}},
      \]
      which depend only on $B$ and $\alpha$ and which have at most $10$ 
      nonzero entries such that
      \[
  z_{s}=y_{s}+W_{2}  \cdot \sigma \left(
  W_{1} \cdot y_{s} + b_{1}
  \right) + b_{2}
  \]
  results in
  \[
  z_s^{(i)}
  =
  \begin{cases}
    y_s^{(i)} & \mbox{ if }
    i \in \{1, \dots, d+l+1\}, \\
    \alpha \cdot ( y_s^{(d+l+4)}-B) & \mbox{ if }
    i = d+l+3, \\
    0 & \mbox{ if } i \in \{d+l+2,d+l+4\}.
  \end{cases}
  \]
      \end{lemma}

    \noindent
        {\bf Proof.} W.l.o.g. we assume $d_{ff}=8$.
        We choose $b_1=0$, $b_2=0$,
        \[
        W_1=
        \left(
    \begin{array}{cccccccccc}
      0 & \dots & 0 & 0 &  0 &\dots & 0 & 1 &  0 & 0 \\
      0 & \dots & 0 & 0 &  0 &\dots & 0 & -1 &  0 & 0 \\
      0 & \dots & 0 & 0 &  0 &\dots & 0 & 0 &  1 & 0 \\
      0 & \dots & 0 & 0 &  0 &\dots & 0 & 0 &  -1 & 0 \\
      0 & \dots & 0 & -B &  0 &\dots & 0 & 0 &  0 & 1 \\
      0 & \dots & 0 & B &  0 &\dots & 0 & 0 &   0 & -1 \\
      0 & \dots & 0 & 0 &  0 &\dots & 0 & 0 &  0 & 1 \\
      0 & \dots & 0 & 0 &  0 &\dots & 0 & 0 &   0 & -1 \\
      \end{array}
    \right),
        \]
where  all columns
except columns number $d+1$, $d+l+2$, $d+l+3$ and $d+l+4$ are zero,
and
        \[
        W_2=
        \left(
    \begin{array}{cccccccc}
      0 & 0 & 0 & 0 & 0 & 0 & 0 & 0\\
      \vdots & &&&&& & \vdots \\
      0 & 0 & 0 & 0 & 0 & 0 & 0 & 0\\
      -1 & 1 & 0 & 0 & 0 & 0 & 0 & 0\\
      0 & 0 & -1 & 1 & \alpha & -\alpha & 0 & 0\\
      0 & 0 & 0 & 0 & 0 & 0 & -1 & 1
      \end{array}
    \right),
        \]
where all rows except row number $d+l+2$, $d+l+3$ and $d+l+4$ are zero.
Then we have
\begin{eqnarray*}
          &&
W_{2} \cdot \sigma \left(
  W_{1} \cdot y_{s} + b_{1}
  \right) + b_{2}
  \\
  &&
  =
  \left(
  \begin{array}{c}
    0 \\
    \vdots \\
    0 \\
    -   ( \sigma(y_{s}^{(d+l+2)})-    \sigma(-y_{s}^{(d+l+2)})) \\
    -   ( \sigma(y_{s}^{(d+l+3)})-    \sigma(-y_{s}^{(d+l+3)})) +
    \alpha \cdot \sigma(y_{s}^{(d+l+4)}-B)
    - \alpha \cdot \sigma(B-y_{s}^{(d+l+4)})\\
    -   ( \sigma(y_{s}^{(d+l+4)})-    \sigma(-y_{s}^{(d+l+4)})) 
    \end{array}
  \right).
        \end{eqnarray*}
        Because of
        \[
\sigma(u)-\sigma(-u)=u
        \]
        for $u \in \R$ this implies the assertion.
        \hfill $\Box$

            \begin{lemma}
      \label{le3b}
      Let $d_{model}=d+l+4$ and
      let $y=(y_1, \dots, y_l)$ for some $y_i \in \R^{d_{model}}$.
      Let $d_{ff} \geq 3$ and let $\alpha \in \R$.
      Then there exists matrices and vectors
      \[
      W_{1} \in \R^{d_{ff} \times d_{model}}, b_{1} \in \R^{d_{ff}},
      W_{2} \in \R^{d_{model} \times d_{ff}}, b_{2} \in \R^{d_{model}},
      \]
       which have at most three 
      nonzero entries such that
      \[
  z_{s}=y_{s}+W_{2}  \cdot \sigma \left(
  W_{1} \cdot y_{s} + b_{1}
  \right) + b_{2}
  \]
  results in
  \[
  z_s^{(i)}
  =
  \begin{cases}
    y_s^{(i)} & \mbox{ if }
    i \in \{1, \dots, d+l+4\} \setminus \{ d+l+2 \}, \\
    \max\{ y_s^{(d+l+4)},0\} & \mbox{ if }
    i = d+l+2 
  \end{cases}
  \]
  $(s \in \{1, \dots, l\})$.
      \end{lemma}

    \noindent
        {\bf Proof.} Follows as in the proof of Lemma \ref{le3}
        by choosing $W_2$, $b_2$, $W_1$ and $b_1$ such that
        \[
W_{2} \cdot \sigma \left(
  W_{1} \cdot y_{s} + b_{1}
  \right) + b_{2}
  =
  \left(
  \begin{array}{c}
    0 \\
    \vdots \\
    0 \\
    -   ( \sigma(y_{s}^{(d+l+2)})-    \sigma(-y_{s}^{(d+l+2)})) +
    \sigma(y_{s}^{(d+l+4)})\\
    0  \\
    0 \\
    \end{array}
  \right)
  \]
  holds. \hfill $\Box$ \\

  \noindent
      {\bf Remark 5.} Assume $d_{ff} \geq 4$.
      It follows from the proof of Lemma \ref{le3b}
      that we can modify $W_2$ such that at most four
      of its entries are nonzero and
      $z_s^{(d+l+2)}=y_s^{(d+l+4)}$ holds instead of
      $z_s^{(d+l+2)}=\max\{ y_s^{(d+l+4)}, 0 \}$.\\

        Next we combine the previous results to construct
        a Transformer encoder which uses $ 2 \cdot M \cdot d$
        pairs 
        of single-attention layers and piecewise feedforward
        layers in order to compute a basis function
        of a tensor product spline space of degree $M$.

        \begin{lemma}
          \label{le4}
          Let $K \in \N$, let $u_k \in \R$ $(k=1, \dots, K-1)$, set
          \[
B_j(x)=x^j \quad \mbox{for } j=0,1,\dots, M
\]
and set
\[
B_j(x)=(x-u_{j-M})_+^M  \quad \mbox{for } j=M+1,M+2,\dots, M+K-1.
\]
Let $j_1, \dots, j_d \in \{0,1,\dots, M+K-1\}$ and $\alpha \in \R$.
Then there exists a transformer network consisting of
$2 \cdot  M \cdot d$ pairs of layers layers,
where in each pair
the first layer is a single-attention layer and
the second layer is a piecewise feedforward
        neural network, where all matrices and vectors
        have at most $10$  nonzero entries and where all matrices
        and vectors depend only on $(u_k)_k$, $B$ and $\alpha$, which gets
        as input $z$ as in Lemma \ref{le1} (with $b_j=1$ $(j=1, \dots, d)$)
        and produces
        as output $z_{ M \cdot d}$ which satisfies
        \[
z_{M \cdot d}^{(d+l+3)} = \alpha \prod_{k=1}^d B_{j_k}(x^{(k)}).
\]
          \end{lemma}

        \noindent
            {\bf Proof.}
            Each $B_j(x)$ can be written as
            \[
B_j(x)= \prod_{k=1}^M B_{j,k}(x)
            \]
where $B_{j,k}(x)$ is one of the functions
            \[
            x \mapsto 1, \quad x \mapsto x \quad
            \mbox{and} \quad x \mapsto (x -u_r)_+. 
            \]            
            Using Lemma \ref{le1},
            Lemma \ref{le2},
            Remark 4, Lemma \ref{le3}, Lemma \ref{le3b}, Remark 5
            and Lemma \ref{le4}
            we can combine two pairs of attention layers and
            piecewise feedforward layers such that they produce from the
            input
            $z_{0,j}$ given as in Lemma \ref{le1} with $b_j \geq 0$
            an output $y_j$ where $y_1^{(d+l+3)}$ is the product
            of $b_j$ and one of the functions
            \[
            x^{(k)} \mapsto 1, \quad x^{(k)} \mapsto x^{(k)} \quad
            \mbox{and} \quad x^{(k)} \mapsto (x^{(k)} -u_r)_+ 
            \]
            and where $y_j^{(s)}$ is equal to $z_{0,j}^{(s)}$ otherwise.
            Using this repeatedly  we get the assertion.
            \hspace*{3cm} \hfill $\Box$

            \subsection{Approximation of piecewise polynomials by transformer networks
  with multi-head attention}
            In this subsection we generalize the results from the
            previous subsection to Transformer encoders with
            multi-head attention. The basic idea is to extend the
            coding of the input by repeating the previous coding $h$ times
            and to define the multi-head attention such that each
            attention unit makes the computation of the previous
            subsection in one of the copies of the original coding.

We use the following coding of the original input:
We represent the input sequence (\ref{se2eq1}) by
\begin{equation}
  \label{se5eq2}
z_0=(z_{0,1}, \dots, z_{0,l}) \in \R^{l \cdot d_{model}},
\end{equation}
where
\begin{equation}
  \label{se5eq3}
  d_{model}
  =
  h \cdot (d+l+4).
\end{equation}
Here we repeat the coding of the previous subsection $h$ times.
More precisely we set for all $k \in \{1, \dots, h\}$
\[
z_{0,j}^{((k-1) \cdot (d+l+4) +s)}
=
\begin{cases}
  x_j^{(s)} & \mbox{if } s \in \{1, \dots, d\} \\
  1 & \mbox{if } s=d+1 \\
  \delta_{s-d-1,j} & \mbox{if } s \in \{d+2, \dots, d+1+l\}\\
  1 & \mbox{if } s=d+l+3 \\
  0 & \mbox{if } s \in \{d+l+2, d+l+4 \}\\
  \end{cases}
\]

            In our next  result we
            compute $h$ basic functions of the truncated power basis
            in parallel by using a multi-head attention with $h$
            attention heads.

            \begin{lemma}
              \label{le5}
              Let $K$, $u_k$, $M$ and $B_j$ be as in Lemma \ref{le4}.
              Let $d_k \geq 2$ and $d_{ff} \geq 10$.
              Let $h \in \N$ and for $s \in \{1, \dots, h\}$ let
              $j_{s,1}, \dots, j_{s,d} \in \{0,1, \dots, M+K-1\}$
              and $\alpha_s \in \R$.
              Then there exists a Transformer encoder consisting of
              $ 2 \cdot M \cdot d$ pairs of layers, where the first layer is a
              multi-head attention layer with $h$ attention units and
              the second layer is a piecewise feedforward
        neural network, and where all matrices and vectors
        have at most $10$  nonzero entries and where all matrices
        and vectors depend only on $(u_k)_k$, $B$ and
        $\alpha_s$ $(s \in \{1, \dots, k\})$, which gets
        as input $z_0$ defined in (\ref{se5eq2}) and produces
        as output $z_{M \cdot d}$ which satisfies
        \[
z_{M \cdot d,1}^{(s-1) \cdot (d+l+4) + (d+l+3)} = \alpha_s \prod_{k=1}^d B_{j_{s,k}}(x^{(k)})
\]
for all $s \in \{1, \dots, h\}$.
              \end{lemma}

            \noindent
                {\bf Proof.}
                The result is a straightforward extension of the
                proof of Lemma \ref{le4}. The basic idea is as follows.
                Each attention head of the network works only on one of the
                $h$ copies of the coding of the previous subsection.
                In each attention unit it makes
                the same computations as in the proof of Lemma \ref{le4},
                using only its special part of the coding of the input.
                \hfill $\Box$

            \begin{lemma}
              \label{le5b}
              Let $K$, $u_k$, $M$ and $B_j$ be as in Lemma \ref{le4}.
              Let $h \in \N$ and for $s \in \{1, \dots, h\}$ let
              $j_{s,1}, \dots, j_{s,d} \in \{0,1, \dots, M+K-1\}$
              and $\alpha_s \in \R$.
              Let $d_k \geq 2$ and $d_{ff} \geq \max\{10, 2 \cdot h + 2 \}$.
              Then there exists a Transformer encoder consisting of
               $ 2 \cdot M \cdot d+1$ pairs of layers, where the first layer is a
              multi-head attention layer with $h$ attention units and
              the second layer is a piecewise feedforward
              neural network, and where all matrices and vectors
              in the first $ 2 \cdot M \cdot d$ pairs of layers
              have at most $10$  nonzero entries,
              where the matrices and vectors in the 
              $ 2 \cdot M \cdot d+1$-th pair of layers have
              together at most $2 \cdot h+2$ 
              nonzero components
              and where all matrices
        and vectors depend only on $(u_k)_k$, $B$ and
        $\alpha_s$ $(s \in \{1, \dots, k\})$, which gets
        as input $z_0$ defined in (\ref{se5eq2}) and produces
        as output $z_{M \cdot d+1}$ which satisfies
        \[
        z_{M \cdot d+1,1}^{h \cdot (d+l+4) } =
        \sum_{s=1}^h
        \alpha_s \prod_{k=1}^d B_{j_{s,k}}(x^{(k)}).
\]
              \end{lemma}

            \noindent
                {\bf Proof.}
                We use the construction of the proof of Lemma \ref{le5} to define
                the first $2 \cdot M \cdot d$ pairs of layers. After that we choose
                $W_{V,M \cdot d}=0$ (which results in $y_{M \cdot d + 1}=z_{M \cdot d}$)
                and choose $W_{1, M \cdot d+1}$, $b_{1, M \cdot d+1}$,
                $W_{2, M \cdot d+1}$, $b_{2, M \cdot d+1}$, 
                such that
                \begin{eqnarray*}
                  &&
                  \hspace*{-0.6cm}
W_{2, M \cdot d+1} \cdot \sigma \left(
  W_{1, M \cdot d+1} \cdot y_{s} + b_{1, M \cdot d+1}
  \right) + b_{2, M \cdot d+1}
  \\[0.3cm]
  &&
                  \hspace*{-0.6cm}
  =
  \left(
  \begin{array}{c}
    0 \\
    \vdots \\
    0 \\
    -   ( \sigma(y_{M \cdot d+1, 1}^{h \cdot (d+l+4) })-    \sigma(-y_{M \cdot d+1,1}^{h \cdot (d+l+h) })) +
    \sigma(
\sum_{s=1}^h y_{M \cdot d+1, 1}^{s \cdot (d+l+4)  }
)
-
   \sigma(-
\sum_{s=1}^h y_{M \cdot d+1, 1}^{s \cdot (d+l+4)  }
)
\\
    \end{array}
  \right)
  \end{eqnarray*}
  holds.
                \hfill $\Box$ \\

                \noindent
                    {\bf Remark 6.}
                    The transformer network in Lemma \ref{le5b} has
                    \begin{eqnarray*}
                      &&
                    N \cdot h \cdot
                    (2 \cdot d_k \cdot d_{model} + d_v \cdot d_{model}
                    + 2 \cdot d_{ff} \cdot d_{model} +
                    d_{ff} + d_{model})
                    \\
                    &&
                    =
                    (2 \cdot M \cdot d + 1) \cdot h \cdot (
                    2 \cdot d_k \cdot h \cdot (d+l+4)
                    +
                    (d+l+4) \cdot h \cdot (d+l+4) \\
                    &&
                    \quad
                    +
                    2 \cdot d_{ff} \cdot h \cdot (d+l+4)
                    + d_{ff} +  h \cdot (d+l+4))
                    \\
                    &&
                    \leq
                    235 \cdot M \cdot (\max\{l,d,d_k,d_{ff} \})^3 \cdot h^2
                    \end{eqnarray*}
                    many parameters, and of these parameters
                    at most $144 \cdot M \cdot d \cdot h$ are not equal to zero.
                    Here the positions are fixed where nonzero parameters
                    are allowed to appear.\\

                    Next we show how we can approximate a function
                    which satisfies
                    a hierarchical composition model
                    by a Transformer encoder.
                    In order to formulate this result, we introduce some
                    additional notation. In order to
compute a function $h_1^{(\kappa)} \in \mathcal{H}(\kappa, \mathcal{P})$
                    one has to compute different hierarchical composition models of some level $i$ $(i\in \{1, \dots, \kappa-1\})$. Let $\tilde{N}_i$ denote the number of hierarchical composition models of level $i$, needed to compute $h_1^{(\kappa)}$. 
Let
\begin{align}
\label{h}
h_j^{(i)}: \R^{d} \to \R 
%\quad ((p_j^{(i)}, K_j^{(i)}) \in \mathcal{P})
\end{align}
be the $j$--th hierarchical composition model of some level $i$ ($j \in \{1, \ldots, \tilde{N}_i\}, i \in \{1, \ldots, \kappa\}$),
that applies a $(p_j^{(i)}, C)$--smooth function $g_j^{(i)}: \R^{K_j^{(i)}} \to \R$ with $p_j^{(i)} = q_j^{(i)} + s_j^{(i)}$, $q_j^{(i)} \in \N_0$ and $s_j^{(i)} \in (0,1]$, where $(p_j^{(i)}, K_j^{(i)}) \in \mathcal{P}$.
  With this notation we can describe
 the computation of $h_1^{(\kappa)}(\bold{x})$ recursively as follows:
    \begin{equation}\label{hji}
  h_j^{(i)}(\bold{x}) =  g_{j}^{(i)}\left(h^{(i-1)}_{\sum_{t=1}^{j-1} K_t^{(i)}+1}(\bold{x}), \dots, h^{(i-1)}_{\sum_{t=1}^j K_t^{(i)}}(\bold{x}) \right)
  \end{equation}
for $j \in \{1, \dots, \tilde{N}_i\}$ and $i \in \{2, \dots, \kappa\}$,
  and
    \begin{equation}\label{hj1}
  h_j^{(1)}(\bold{x}) = g_j^{(1)}\left(x^{\left(\pi(\sum_{t=1}^{j-1} K_t^{(1)}+1)\right)}, \dots, x^{\left(\pi(\sum_{t=1}^{j} K_t^{(1)})\right)}\right) 
    \end{equation}
    holds
        for $j \in \{1, \dots, \tilde{N}_1 \}$
  for some function $\pi: \{1, \dots, \tilde{N}_1\} \to \{1, \dots, d\}$. 
  %Then for some $i=l, l-1, \dots, 1$ the recursion 
  Here
  the recursion
\begin{align}
\label{N}
\tilde{N}_l = 1 \ \text{and} \ \tilde{N}_{i} = \sum_{j=1}^{\tilde{N}_{i+1}} K_j^{(i+1)} 
\quad (i \in \{1, \dots, \kappa-1\})
\end{align}
holds.

                    \begin{theorem}
                      \label{th2}
                      Let $m: \mathbb{R}^d \to \mathbb{R}$ be contained in the class
                      $\mathcal{H}(\kappa, \mathcal{P})$ for some $\kappa \in \N$ and $\mathcal{P} \subseteq [1,\infty) \times \N$.  Let $\tilde{N}_i$ be defined as in \eqref{N}. Each $m$ consists of different functions $h_j^{(i)}$ $(j \in \{1, \ldots, \tilde{N}_i\},$
 $ i\in \{1, \dots, \kappa\})$ defined as in \eqref{h}, \eqref{hji} and \eqref{hj1}. 
 Assume that the corresponding functions $g_j^{(i)}$ are Lipschitz continuous with Lipschitz constant $C_{Lip} \geq 1$ and satisfy
  \begin{equation*}
  \|g_j^{(i)}\|_{C^{q_j^{(i)}}(\R^d)} \leq c_{9}
  \end{equation*}
  for some constant $c_{9} >0$. Denote by $K_{max} = \max_{i,j} K_j^{(i)} < \infty$ the maximal input dimension and set $q_{max} = \max_{i,j} q_j^{(i)} < \infty$,
  where $q_j^{(i)}$ is the integer part of the smoothness $p_j^{(i)}$
  of $g_j^{(i)}$.  Let $A \geq 1$.
  Choose $h \in \N$ such that
  \begin{equation}
    \label{th2eq1}
h \geq c_{10}
\end{equation}
holds for some sufficiently large constant $c_{10}$, choose
  \[
I \geq \sum_{i=1}^\kappa \tilde{N}_i \quad \mbox{and} \quad d_k \geq 2
\]
  and set
  \[
  N= I \cdot (2 \cdot q_{max} \cdot K_{max} +1), \quad
  d_{model}=d_v=h \cdot I \cdot (d+l+4)
  \]
  and
  \[
  L_n= 144 \cdot (q_{max}+1) \cdot K_{max} \cdot I \cdot h. 
  \]
  Then there exists a transformer network $f_{\bar{\theta}}$
with $\|\bar{\theta}\|_0 \leq L_n$ which satisfies
  \[
  \|f_{\bar{\theta}}-m\|_{\infty, [-A,A]^d} \leq c_{11} \cdot (K_{max}+1)^\kappa \cdot
  \max_{j,i} h^{-p_j^{(i)}/K_j^{(i)}}.
  \]
                    \end{theorem}

                    \noindent
                        {\bf Proof.}
                        From the Lipschitz continuity of the
                        $g_j^{(i)}$ and the recursive
                        definition of the $h_j^{(i)}$ we can conclude
                        that there exists $\bar{A} \geq A$ such that
                        \begin{equation}
                          \label{pth2eq1}
h_j^{(i)}(x) \in [-\bar{A},\bar{A}]
\end{equation}
holds for all $x \in [-A,A]^d$, $j \in \{1, \dots, \tilde{N}_i\}$
and $i \in \{1, \dots, \kappa-1\}$.

                        Our Transformer encoder successively
                        approximates 
                        $h_1^{(1)}(x)$, \dots, $h_{N_1}^{(1)}(x)$,
                        $h_1^{(2)}(x)$, \dots, $h_{N_2}^{(2)}(x)$,
                        \dots, $h^{(\kappa)}_1(x)$.
                        Here $h_i^{(j)}$ is
                        approximated by computing
                        in a first step truncated power basis
                        of a tensorproduct spline space
                        of degree $q_{i}^{(j)}$ on a equidistant grid in
                        \[
                          [-\bar{A}-1,\bar{A}+1]^{K_i^{(j)}}
                          \]
                          consisting of $h$ basis
                        functions, which are evaluated at the arguments
                        of $h_i^{(j)}$ in (\ref{hj1}),
                        and by using in a second step a linear combination of these basis functions to approximate
                        \[
                        g_{j}^{(i)}\left(h^{(i-1)}_{\sum_{t=1}^{j-1} K_t^{(i)}+1}(\bold{x}), \dots, h^{(i-1)}_{\sum_{t=1}^j K_t^{(i)}}(\bold{x}) \right)
                        .
                        \]
                        The computation of this truncated power basis can be done as in Lemma \ref{le5b} using layers $(\tilde{N}_{i-1}+j-1) \cdot (2 \cdot q_{max} \cdot K_{max} +1)+1$
                        till
$(\tilde{N}_{i-1}+j) \cdot (2 \cdot q_{max} \cdot K_{max} +1)$
                        of our Transformer encoder
                        and proceeding otherwise as in Lemma \ref{le5}.
                        Using standard approximation results from spline
                        theory (cf., e.g., Theorem 15.2 and
                        proof of Theorem 15.1 in Gy\"orfi et al. (2002)
                        and Lemma 1 in Kohler (2014))
                        this results in an approximation
                        \[
\tilde{g}_{j}^{(i)}
                        \]
                        of $g_{j}^{(i)}$ which satisfies
                        \begin{equation}
                          \label{pth2eq2}
                        \| \tilde{g}_{j}^{(i)} - g_{j}^{(i)} \|_{\infty, [-\bar{A}-1,\bar{A}+1]^{K_j^{(i)}}}
                        \leq
                        c_{12}  \cdot h^{-p_{j}^{(i)}/K_j^{(i)}}.
                        \end{equation}
                        The approximation
                        $\tilde{h}_1^{(\kappa)}(\bold{x})$
                        of $h_1^{(\kappa)}(\bold{x})$ which our Transformer encoder computes
                        is defined  as follows:
    \[
  \tilde{h}_j^{(1)}(\bold{x}) = \tilde{g}_j^{(1)}\left(x^{\left(\pi(\sum_{t=1}^{j-1} K_t^{(1)}+1)\right)}, \dots, x^{\left(\pi(\sum_{t=1}^{j} K_t^{(1)})\right)}\right)
  \]
for  $j \in \{1, \dots, \tilde{N}_1\}$
  and
  \[
    \tilde{h}_j^{(i)}(\bold{x}) =  \tilde{g}_{j}^{(i)}\left(
    \tilde{h}^{(i-1)}_{\sum_{t=1}^{j-1} K_t^{(i)}+1}(\bold{x}), \dots,
    \tilde{h}^{(i-1)}_{\sum_{t=1}^j K_t^{(i)}}(\bold{x}) \right)
  \]
  for $j \in \{1, \dots, \tilde{N}_i\}$ and $i \in \{2, \dots, \kappa\}$.

  From (\ref{th2eq1}), (\ref{pth2eq1}) and (\ref{pth2eq2}) we can
  conclude 
  \[
  |\tilde{h}_j^{(i)}(\bold{x})| \leq
  | \tilde{h}_j^{(i)}(\bold{x}) - h_j^{(i)}(\bold{x})|
  + |h_j^{(i)}(\bold{x})| \leq \bar{A}+1.
  \]
  Consequently we get from  (\ref{pth2eq2})
  \begin{eqnarray*}
    &&
    |
    \tilde{h}_j^{(i)}(\bold{x})
    -
    h_j^{(i)}(\bold{x})
    |
    \\
    &&
\leq
|
\tilde{g}_{j}^{(i)}\left(
    \tilde{h}^{(i-1)}_{\sum_{t=1}^{j-1} K_t^{(i)}+1}(\bold{x}), \dots,
    \tilde{h}^{(i-1)}_{\sum_{t=1}^j K_t^{(i)}}(\bold{x}) \right)
    -
g_{j}^{(i)}\left(
    \tilde{h}^{(i-1)}_{\sum_{t=1}^{j-1} K_t^{(i)}+1}(\bold{x}), \dots,
    \tilde{h}^{(i-1)}_{\sum_{t=1}^j K_t^{(i)}}(\bold{x}) \right)
    |
    \\
    &&
    \quad
    +
    |
g_{j}^{(i)}\left(
    \tilde{h}^{(i-1)}_{\sum_{t=1}^{j-1} K_t^{(i)}+1}(\bold{x}), \dots,
    \tilde{h}^{(i-1)}_{\sum_{t=1}^j K_t^{(i)}}(\bold{x}) \right)
    -
    g_{j}^{(i)}\left(
    h^{(i-1)}_{\sum_{t=1}^{j-1} K_t^{(i)}+1}(\bold{x}), \dots,
    h^{(i-1)}_{\sum_{t=1}^j K_t^{(i)}}(\bold{x}) \right)
    |
    \\
    &&
    \leq
                            c_{12} \cdot h^{-p_{j}^{(i)}/K_j^{(i)}}
                            +
                            \\
                            &&
                            \quad
    |
g_{j}^{(i)}\left(
    \tilde{h}^{(i-1)}_{\sum_{t=1}^{j-1} K_t^{(i)}+1}(\bold{x}), \dots,
    \tilde{h}^{(i-1)}_{\sum_{t=1}^j K_t^{(i)}}(\bold{x}) \right)
    -
    g_{j}^{(i)}\left(
    h^{(i-1)}_{\sum_{t=1}^{j-1} K_t^{(i)}+1}(\bold{x}), \dots,
    h^{(i-1)}_{\sum_{t=1}^j K_t^{(i)}}(\bold{x}) \right)
    |
    \\
    &&
    \leq
                            c_{12} \cdot h^{-p_{j}^{(i)}/K_j^{(i)}}
                            + c_{13} \cdot \sum_{s=1}^{K_j^{(i)}}
                            |\tilde{h}^{(i-1)}_{\sum_{t=1}^{j-1} K_t^{(i)}+s}
                            (\bold{x})
                            - h^{(i-1)}_{\sum_{t=1}^{j-1} K_t^{(i)}+s}
                            (\bold{x})
         |,
    \end{eqnarray*}
  where the last inequality holds due to the
  the Lipschitz continuity of $g_{j}^{(i)}$. Together with
  \[
      |
    \tilde{h}_j^{(1)}(\bold{x})
    -
    h_j^{(1)}(\bold{x})
    |
    \leq
     c_{14} \cdot h^{-p_{j}^{(1)}/K_j^{(1)}},
\]
which follows again from (\ref{pth2eq2}), an easy induction shows
\[
|    \tilde{h}_1^{(\kappa)}(\bold{x})
    -
    h_1^{(\kappa)}(\bold{x})
    |
    \leq
    c_{15} \cdot (K_{max}+1)^\kappa \cdot
\max_{j,i} h^{-p_j^{(i)}/K_j^{(i)}} .    
\]

                        \hfill $\Box$
                        
\subsection{A bound on the covering number}

In this subsection we prove the following bound on the covering
number.

\begin{lemma}
  \label{le6}
  Let $\F$ be the set of all functions
  \[
(x_1, \dots, x_l) \mapsto z_N \cdot w + b,
  \]
  where $z_N$ is defined in Section \ref{se2} depending on
  \begin{equation}
    \label{se5eq4}
  (W_{Q,r,s}, W_{K,r,s}, W_{V,r,s})_{r \in \{1, \dots, N\}, s \in \{1, \dots, h \}},
    \quad
    (W_{r,1},b_{r,1},W_{r,2},b_{r,2})_{r \in \{1, \dots, N\}}
  \end{equation}
  and where the total number of nonzero components in (\ref{se5eq4}) and
  $w$ and $b$ is bounded by $L \in \N$. Let $\beta \geq 0$ and let $T_\beta \F$
  be the set of all functions in $\F$ truncated on height $\beta$ and $-\beta$.
  Then we have for any $0<\epsilon<\beta/2$
  \[
  \sup_{z_1^n \in (\R^{d \cdot l})^n)} \log \Nu(\epsilon, T_\beta \F, z_1^n) \leq
  c_{16} \cdot L \cdot N^2 \cdot \log( \max\{N,h,d_{ff},d_k,l\}) \log \left(
\frac{\beta}{\epsilon}
  \right).
  \]
  \end{lemma}

In order to prove Lemma \ref{le6} we will first show the following
bound on the VC-dimension of subsets of $\F$, where the nonzero
components appear only at fixed positions.

\begin{lemma}
  \label{le7}
   Let $\G$ be the set of all functions
  \[
(x_1, \dots, x_l) \mapsto z_N \cdot w + b,
  \]
  where $z_N$ is defined in Section \ref{se2} depending on
   (\ref{se5eq4})
  and where there are at most $L$ fixed components in (\ref{se5eq4}) and
  $w$ and $b$ where the entries are allowed to be nonzero. Then we have
  \[
V_{\F^+} \leq c_{17} \cdot L \cdot N^2 \cdot \log( \max\{N,h,d_{ff},l\})
  \]
  \end{lemma}

The proof of Lemma \ref{le7} is a modification of the proof of
Theorem 6 in Bartlett et al. (1999).
In the proof of Lemma \ref{le7} we will need the following two auxiliary
results.

\begin{lemma}
\label{le8}
Suppose $W\leq m$ and let $f_1,...,f_m$ be polynomials of degree at most $D$ in $W$ variables. Define
\[
K\coloneqq|\{\left(\sgn(f_1(\ba)),\dots,\sgn(f_m(\ba))\right) : \ba\in\R^{W}\}|.
\]
Then we have 
\[
K\leq2\cdot\left(\frac{2\cdot e\cdot m\cdot D}{W}\right)^{W}.
\]
\end{lemma}
\noindent
    {\bf Proof.} See Theorem 8.3 in Anthony and Bartlett (1999).
		\hfill $\Box$

      \begin{lemma}
\label{le9}
	Suppose that $2^m\leq2^L\cdot(m\cdot R/w)^w$ for some $R\geq16$ and $m\geq w\geq L\geq0$. Then,
	\[
	m\leq L+w\cdot\log_2(2\cdot R\cdot\log_2(R)).
	\]
\end{lemma}
\noindent
{\bf Proof.}
See Lemma 16 in Bartlett et al. (2019).
\hfill $\Box$          

\noindent
    {\bf Proof of Lemma \ref{le7}}.
    Let $\HH$ be the set of all functions $h$ defined by
    \[
    h:(\R^d)^l \times \R \rightarrow \R,
    \quad
    h(x,y)=g(x)-y
    \]
    for some $g \in \G$.
                    Let $(x_1, y_1)$, \dots
                    $(x_m,y_m) \in \R^{d \cdot l} \times \R$ be such that
                    \begin{equation}
                      \label{ple7eq1}
                      |
                      \{
(sgn(h(x_1,y_1)), \dots, sgn(h(x_m,y_m))) \, : \, h \in \HH
                      \}
                      |
                      =2^m.
                      \end{equation}
It suffices to show
                         \begin{equation}
                           \label{ple7eq2}
                           m \leq c_{17} \cdot L \cdot N^2 \cdot \
                           \log( \max\{N,h,d_{ff},l\}).
                      \end{equation}
                         To show this we partition $\G$ in subsets
                         such that for each subset all
                         \[
g(x_i) \quad (i=1, \dots, m)
                         \]
                         are polynomials of some fixed degree and use Lemma \ref{le8}
                         in order to derive an upper bound on the left-hand side of (\ref{ple7eq1}).
                         This upper bound will depend polynomially on $m$ which will enable
                         us to conclude (\ref{ple7eq2}) by an application of Lemma \ref{le9}.
                         
                         Let
                         \[
                         \bar{\theta} = \left(
  (W_{Q,r,s}, W_{K,r,s}, W_{V,r,s})_{r \in \{1, \dots, N\}, s \in \{1, \dots, h\}},
    (W_{r,1},b_{r,1},W_{r,2},b_{r,2})_{r \in \{1, \dots, N\}}, w, b
                         \right)
                         \]
                         be the parameters which determine a function in
                         $\F$. By assumption, each function in $\G$
                         can be also described by such a parameter vector.
                         Here
                         only $L$ components of the matrices and vectors
                         occuring in the parameter vector
                         are allowed to be nonzero
                         and the positions
                         where these nonzero parameters can occur are fixed.
                         Denote the vector in $\R^L$ which contains all
                         values of these possible nonzero parameters
                         by $\theta$. Then we can write
                         \[
                         \G = \{ g(\cdot, \theta):
                         \R^{d \cdot l} \rightarrow \R \, : \,
                         \theta \in \R^L
                         \}.
                         \]
                         In the sequel we construct a partition $\P_{N+1}$ of $\R^L$
                         such that for all $S \in \P_{N+1}$ we have
                         that
                         \[
                         g(x_1, \theta), \dots, g(x_m,\theta)
                         \]
                          (considered as functions of $\theta$) are 
                         polynomials of degree at most $8^{N+1}$ for $\theta \in S$.

                         In order to construct this partition we construct
                         recursively partionions $\P_{0}$, \dots, $\P_N$ of $\R^L$ such that for each $r \in \{1, \dots, N\}$ and all $S \in \P_r$
                         all components in
                         \[
z_r
                         \]
                         (considered as a function of $\theta$) are polynomials of degree at most $8^r$ in $\theta$ for
                         $\theta \in S$.

                         Since all components of $z_0$ are constant as
                         functions of $\theta$ this holds if we set $\P_{0}=\{ \R^L\}$.

                         Let $r \in \{1, \dots, N\}$ and assume that
                         for 
 all $S \in \P_{r-1}$
                         all components in
                         \[
z_{r-1}
                         \]
                         (considered as a function of $\theta$) are polynomials of degree at most $8^{r-1}$ in $\theta$ for
                         $\theta \in S$.
                         Then all components in
                         \[
q_{r-1,s,i}, k_{r-1,s,i} \quad \mbox{and} \quad v_{r-1,s,i}
                         \]
                         are on each set $S \in \P_{r-1}$ polynomials of degree $8^{r-1} +1$.
                         Consequently, for each $S \in \P_{r-1}$ each value
                         \[
<q_{r-1,s,i}, k_{r-1,s,j}> 
                         \]
                         is (considered as a function of $\theta$)
                         a polynomial of degree at most $2 \cdot 8^{r-1}+2$ 
                          for
                          $\theta \in S$. Application of
                          Lemma \ref{le8} yields that
                          \[
                          <q_{r-1,s,i}, k_{r-1,s,j_1}>
                          -
                          <q_{r-1,s,i}, k_{r-1,s,j_2}>
                          \quad
                          (s \in \{1, \dots, h\}, i, j_1, j_2 \in \{1, \dots, l\})
                          \]
                          has at most
                          \[
                          \Delta
                          =
                          2 \cdot
                          \left(
                          \frac{
2 \cdot e \cdot h \cdot l^3 \cdot (2 \cdot 8^{r-1}+2)
                          }{L}
                          \right)^L
                          \]
                          difference sign patterns. If we partiton in each set
                          in $\P_{r-1}$ according to these sign patterns
                          in $\Delta$ subsets, then on each set in the new partition
                          all components in
                          \[
v_{r-1,s,\hat{j}_{r-1,s,i}} \cdot  <q_{r-1,s,i}, k_{r-1,s,\hat{j}_{r-1,s,i}}>
                          \]
                          are polynomials of degree at most
                          $3 \cdot 8^{r-1} + 3$
                          (since on each such set $<q_{r-1,s,i}, k_{r-1,s,\hat{j}_{r-1,s,i}}>$ is equal to one of the $<q_{r-1,s,i}, k_{r-1,s,j}>$).
                          On each set within this partition every
                          component of the $\R^{d_{ff}}$-valued vectors
                          \[
                          W_{r,1} \cdot y_{r,s} + b_{r,1}
                          \quad (s=1, \dots, h)
                          \]
                          is (considered as a function of $\theta$)
                          a polynomial of degree at most $3 \cdot 8^{r-1} +4$.

                          By another application of Lemma \ref{le8}
                          we can refine each set in this partition
                          into
                          \[
2 \cdot
                          \left(
                          \frac{
2 \cdot e \cdot h \cdot d_{ff} \cdot (3 \cdot 8^{r-1} + 4)
                          }{L}
                          \right)^L
                          \]
                          sets such that all components in
                          \begin{equation}
                            \label{ple7eq3}
W_{r,1} \cdot y_{r,s} + b_{r,1} 
                          \end{equation}
                          have the same sign patterns within the refined partition.
                          We call this refined partion $\P_r$.
                          Since on each set of $\P_r$ the sign of all components
                          in (\ref{ple7eq3}) does not change we can conclude
                          that
 all components in
                          \begin{equation}
                            \label{ple7eq3b}
\sigma( W_{r,1} \cdot y_{r,s} + b_{r,1})
                          \end{equation}                         
                          are either equal to zero or they are equal
                          to a polynomial of degree at most
                          $ 3 \cdot 8^{r-1} + 4$.
                          Consequently we have that on each set
                          in $\P_r$ all components of
                          \[
z_r
                          \]
                          are equal to a polynomial of degree at most
$ 3 \cdot 8^{r-1} + 5 \leq 8^r$.

                          Using $\P_{N+1}=\P_N$ we have constructed
                          a partition with
                          \[
                          |\P_{N+1}| = \prod_{r=1}^N
                          \frac{|\P_r|}{|\P_{r-1}|}
                          \leq
                          \prod_{r=1}^N
                          2 \cdot
                          \left(
                          \frac{
2 \cdot e \cdot h \cdot l^3 \cdot 8^r
                          }{L}
                          \right)^L
                          \cdot
2 \cdot
                          \left(
                          \frac{
2 \cdot e \cdot h \cdot d_{ff} \cdot 8^r
                          }{L}
                          \right)^L
                          \]
                          such that for each set in this partition
                          for all $(x,y) \in \{(x_1,y_1), \dots, (x_m,y_m)\}$
\[
g(x)=z_N \cdot w + b \quad \mbox{and} \quad h(x,y)=z_N \cdot w + b -y 
\]
                         (considered as a function of $\theta$) are polynomials of degree at most $8^N+1 \leq 8^{N+1}$ in $\theta$ for
                         $\theta \in S$.

Using
\begin{eqnarray*}
&&
|
                      \{
(sgn(h(x_1,y_1)), \dots, sgn(h(x_m,y_m))) \, : \, h \in \HH
                      \}
                      |
                      \\
                      &&
                      \leq 
                      \sum_{S \in \P_{N+1}}
                      |
                      \{
(sgn(g(x_1, \theta)-y_1), \dots, sgn(g(x_m, \theta)-y_m)) \, : \, \theta \in S
                      \}
                      |
                      \end{eqnarray*}
                      we can apply one more time Lemma \ref{le8} to conclude
                      \begin{eqnarray*}
                        2^m &=&
  |
                      \{
(sgn(h(x_1,y_1)), \dots, sgn(h(x_m,y_m))) \, : \, h \in \HH
                      \}
                      |
                      \\
                      &
                      \leq&
                      |\P_{N+1}| \cdot
                      2 \cdot \left(
\frac{2 \cdot e \cdot m \cdot 8^{N+1}}{L}
\right)^L
\\
&
\leq
&                     
                      2 \cdot \left(
\frac{2 \cdot e \cdot m \cdot 8^{N+1}}{L}
\right)^L
\cdot
\prod_{r=1}^N
                          2 \cdot
                          \left(
                          \frac{
2 \cdot e \cdot h \cdot l^3 \cdot 8^r
                          }{L}
                          \right)^L
                          \cdot
2 \cdot
                          \left(
                          \frac{
2 \cdot e \cdot h \cdot d_{ff} \cdot 8^r
                          }{L}
                          \right)^L
                          \\
                          &
                          \leq &
                          2^{2 \cdot N+1}
                          \cdot
                          \left(
                          \frac{
m \cdot 2 \cdot e \cdot (2 N +1)  \cdot  h \cdot (\max\{l, d_{ff}\})^3 \cdot 8^{N+1}
                          }{(2 N +1) \cdot L}
                          \right)^{(2 N +1) \cdot L}.
                        \end{eqnarray*}
                      Assume $m \geq  (2 N +1) \cdot L $. Application of Lemma \ref{le9} with $L= 2 \cdot N +1$,
                      $R= 2 \cdot e \cdot (2 N +1)  \cdot  h \cdot (\max\{l,d_{ff}\})^3 \cdot 8^{N+1} \geq 16$ and
                      $w=(2 N +1) \cdot L$ yields
                        \[
                        m \leq (2 \cdot N +1)
                        + (2 \cdot N +1) \cdot L \cdot\log_2(2\cdot R\cdot\log_2(R)) \leq c_{18} \cdot L \cdot N^2 \cdot \log( \max\{N,h,d_{ff},l\}),
                        \]
                     which implies (\ref{ple7eq2}).
         	    \hfill $\Box$

                    \noindent
                        {\bf Proof of Lemma \ref{le6}.}
                        The functions in the function set $\F$ depend on
                        at most
                        \[
\lceil c_{19} \cdot N \cdot h^2 \cdot (\max\{d_k,d_{ff},d,l\})^3 \rceil
                        \]
                        many parameters, and of these parameters at most $L$ are allowed
                        to be nonzero. We have
                        \[
\left( {\lceil c_{19} \cdot N \cdot h^2 \cdot (\max\{d_k,d_{ff},d,l\})^3 \rceil
  \atop L} \right)
\leq
\left( \lceil c_{19} \cdot N \cdot h^2 \cdot (\max\{d_k,d_{ff},d,l\})^3 \rceil
\right)^L
\]
many possibilities to choose these positions. If we fix these
positions, we get one function space $\G$ for which we can bound its
VC dimension by Lemma \ref{le7}.
                    Using Lemma \ref{le7},
$V_{T_{\beta} \G^+}
\leq
V_{\G^+},$
Lemma 9.2 and Theorem 9.4
in Gy\"orfi et al. (2002) we get
\begin{eqnarray*}
\mathcal{N}_1 \left(\epsilon,   T_{\beta} \G,
   \bx_1^n\right)
&
   \leq
   &
3 \cdot \left(
\frac{4 e \cdot \beta}{\epsilon}
\cdot
\log
\frac{6 e \cdot \beta}{\epsilon}
\right)^{V_{T_{\beta} \G^+}}
\\
&
\leq &
3 \cdot \left(
\frac{6 e \cdot \beta}{\epsilon}
\right)^{
c_{20} \cdot L \cdot N^2 \cdot \log( \max\{N,h,d_{ff},l\})
}.
\end{eqnarray*}
From this we can conclude
\begin{eqnarray*}
  &&
  \sup_{z_1^n \in (\R^{d \cdot l})^n)} \log \Nu(\epsilon, T_\beta \F, z_1^n)
  \\
  &&
  \leq
  L \cdot
  \log
\left( \lceil c_{19} \cdot N \cdot h^2 \cdot (\max\{d_k,d_{ff},d,l\})^3 \rceil
 \right)
\\
&& \quad
 +
c_{20} \cdot L \cdot N^2 \cdot \log( \max\{N,h,d_{ff},l\})
\cdot \log \left(
\frac{\beta}{\epsilon}
\right)
\\
&&
\leq
c_{21} \cdot L \cdot N^2 \cdot \log( \max\{N,h,d_{ff},d_k,l\})
\cdot \log \left(
\frac{\beta}{\epsilon}
\right).
  \end{eqnarray*}
\hfill $\Box$


\begin{thebibliography}{*}

 \bibitem{AnBa99}
Anthony, M., and Bartlett, P. L. (1999).
 {\it Neural Network Learning: Theoretical Foundations}.
 Cambridge University Press.


\bibitem{BK17}
  Bauer, B., and Kohler, M. (2019).
On deep learning as a remedy for the curse of dimensionality in nonparametric regression. {\it Annals of Statistics}, {\bf 47}, pp. 2261-2285.

 
\bibitem{BHLM19}
  Bartlett, P. L., Harvey, N., Liaw, C., and Mehrabian, A. (2019).
  Nearly-tight VC-dimension bounds for piecewise linear neural networks.
  {\it Journal of Machine Learning Research}, {\bf 20}, pp. 1--17.


\bibitem{Co68}
Cover, T. M. (1968).
Rates of convergence of nearest neighbor procedures.
In {\it Proceedings of the Hawaii International Conference on
Systems Sciences}, pp. 413-415, Honolulu, HI.


\bibitem{DCLT19} 
Devlin, J., Chang, M.-W., Lee, K., and
Toutanova, K. (2019).
BERT: Pre-training of Deep Bidirectional Transformers for
Language Understanding.
arXiv: 1810.04805.


\bibitem{Dev82}
Devroye, L. (1982).
Necessary and sufficient conditions for the almost everywhere
convergence of nearest neighbor regression function estimates.
{\it Zeitschrift f\"ur Wahrscheinlichkeitstheorie und verwandte
Gebiete}, {\bf 61}, pp. 467-481.


\bibitem{DGL96}
Devroye, L., Gy\"orfi, L., and Lugosi, G. (1996).
{\it A Probabilistic Theory Of Pattern Recognition}.
Springer.


  
 \bibitem{GKKW02}
 Gy\"orfi, L., Kohler, M., Krzy\.zak, A., and Walk, H.\ (2002).
 {\it A Distribution-Free Theory of Nonparametric Regression}.
 Springer.

 

\bibitem{ImFu18}
Imaizumi, M.\ and Fukamizu, K.\ (2019).
Deep neural networks learn non-smooth functions effectively.
{\it Proceedings of the 22nd International Conference on Artificial 
Intelligence and Statistics (AISTATS 2019)}.
Naha, Okinawa, Japan.

\bibitem{Ko14}
  Kohler, M. (2014).
  Optimal global rates of convergence for noiseless regression estimation problems with adaptively chosen design.
  {\ it Journal of Multivariate Analysis}, {\bf 132}, pp. 197-208. 

\bibitem{KoKr15}
Kohler, M.\ and Krzy\.zak, A.\ (2017).
Nonparametric regression based on hierarchical interaction models.
\textit{IEEE Transaction on Information Theory}, \textbf{63}, 
pp.\ 1620--1630.

\bibitem{KoLa18}
Kohler, M.\ and Langer, S.\ (2020).
On the rate of convergence of fully connected deep neural network 
regression estimates. To appear in {\it Ann. Stat.}
arXiv: 1908.11133


\bibitem{La20}
Langer, S.\ (2021a).
Analysis of the rate of convergence of fully connected deep neural 
network
regression estimates with smooth activation function, 
{\it Journal of Multivariate Analysis}, {\bf 182}, pp. 104695

\bibitem{La20b}
Langer, S.\ (2021b).
Approximating smooth functions by deep neural networks with 
sigmoid activation function, {\it Journal of Multivariate Analysis}, 
{\bf 182}, pp.  104696


\bibitem{LSYZ21}
Lu, J., Shen, Z., Yang, H. and Zhang, S. (2020)
Deep Network Approximation for Smooth Functions.
arxiv: 2001.03040


\bibitem{Sch19}
Schmidt-Hieber, J.\ (2019).
Deep ReLU networks approximation of functions on a manifold.
arxiv:1908.00695


\bibitem{BERT19}
Sun, C., Qiu, X., Xu, Y., and Huang, X. (2019). 
How to fine-tune BERT for text classification? {\it In China National Conference on Chinese Computational Linguistics}, pp. 194-206. Springer, Cham.


\bibitem{S18}
Suzuki, T.\ (2018).
Adaptivity of deep ReLU network for learning in Besov and mixed 
smooth Besov spaces: optimal rate and curse of dimensionality.
arXiv: 1810.08033.

\bibitem{SN19}
Suzuki, T.\ and Nitanda, A.\ (2019).
Deep learning is adaptive to intrinsic dimensionality of
model smoothness in anisotropic Besov space.
arXiv: 1910.12799.

\bibitem{Vasw17}
 Vaswani, A., Shazeer, N., Parmar, N., Uszkoreit, J., Jones, L., Gomez, A.,  Kaiser, L., and  Polosukhin,I. (2017).
 Attention is all you need.
Arxiv 1706.03762.


\bibitem{Ya18}
Yarotsky, D.\ (2017).
Error bounds for approximations with deep ReLU networks, 
{\it Neural Networks},  \textbf{94}, pp.\ 103--114.

\bibitem{YaZh19}
Yarotsky, D.\ and Zhevnerchuk, A.\ (2020).
The phase diagram of approximation rates for deep
neural networks.
In {\it Advances in Neural Information Processing Systems}, 
\textbf{33}, 
pp.\ 13005--13015.


\bibitem{YHPC18}
Young, T., Hazarika, D., Poria, S., and Cambria, E. (2018).
Recent trends in deep learning based natural language processing.
{\it IEEE Computational Intelligence Magazine}, {\bf 13},
pp. 55-75.



\bibitem{CNN15}
Zhang, X., Zhao, J., and LeCun, Y. (2015). 
Character-level convolutional networks for text classification. {\it Advances in neural information processing systems}, {\bf 28}, 649-657.


\end{thebibliography}
\end{document}